\documentclass[11pt,a4paper]{article}
\usepackage{amsmath,amssymb}

\usepackage{amssymb}
\usepackage{amsmath}

\usepackage{color}
\usepackage[colorlinks,linkcolor=blue,citecolor=red]{hyperref}

\DeclareMathSymbol{\bbbr}{\mathalpha}{AMSb}{"52}
\DeclareMathSymbol{\bbbc}{\mathalpha}{AMSb}{"52}

\newcommand\com[1]{}
\newcommand{\slg}{\mathfrak{sl}}
\newcommand{\gl}{\mathfrak{gl}}

\newcommand{\tr}{\operatorname{tr}}
\newcommand\op[1]{\mathop{\rm #1}\nolimits}

\newtheorem{theorem}{Theorem}

\begin{document}

\title{On  integrability in Grassmann geometries: integrable systems associated with fourfolds  in ${\bf Gr}(3, 5)$}

\author{B. Doubrov$^1$, E.V. Ferapontov$^2$, B. Kruglikov$^3$, V.S.  Novikov$^2$}
     \date{}
     \maketitle
     \vspace{-5mm}
\begin{center}
$^1$Department of Mathematical Physics\\
Faculty of Applied Mathematics\\
 Belarussian State University\\
  Nezavisimosti av. 4, 220030 Minsk, Belarus\\
  \ \\
$^2$Department of Mathematical Sciences \\ Loughborough University \\
Loughborough, Leicestershire LE11 3TU \\ United Kingdom \\
\ \\
$^3$Institute of Mathematics and Statistics\\
NT-Faculty \\
University of Troms\o\\
Troms\o\ 90-37, Norway\\
\ \\
e-mails: \\[1ex]  \texttt{doubrov@islc.org} \\
\texttt{E.V.Ferapontov@lboro.ac.uk}\\
\texttt{boris.kruglikov@uit.no} \\
\texttt{V.Novikov@lboro.ac.uk}
\end{center}

\newpage

\begin{abstract}

Let ${\bf Gr}(d, n)$ be the Grassmannian of $d$-dimensional linear subspaces of an $n$-dimensional vector space $V^n$. A submanifold $X\subset {\bf Gr}(d, n)$ gives rise to a differential system $\Sigma(X)$ that governs $d$-dimensional submanifolds of $V^n$ whose Gaussian image is contained in  $X$. Systems of the form  $\Sigma(X)$
appear in numerous applications in continuum mechanics, theory of integrable systems, general relativity and differential geometry, and include such well-known examples as the dispersionless Kadomtsev-Petviashvili equation, the Boyer-Finley equation, Pleba\'nsky's heavenly equations, and so on. 

In this paper we concentrate on the particularly interesting case of this construction which however contains all essential challenges of the general theory, namely the case when $X$ is a fourfold in ${\bf Gr}(3, 5)$. Our main goal is to investigate differential-geometric and integrability aspects of the corresponding systems $\Sigma(X)$. There exist several  approaches to dispersionless integrability:

\begin{itemize}

\item The method of hydrodynamic reductions, based on a decomposition of a given multi-dimensional system into a collection of commuting $(1+1)$-dimensional systems of hydrodynamic type. Exact  solutions defined by these reductions
can be considered as natural dispersionless
analogues of finite-gap solutions of integrable soliton equations. 
\item The method of dispersionless Lax pairs, that is,  pairs of Hamilton-Jacobi type equations whose compatibility conditions 
are equivalent to the system $\Sigma(X)$.
\item Integrability {\it on solutions}, based on the requirement that the characteristic variety of an integrable system  $\Sigma(X)$ defines an integrable background geometry  (Einstein-Weyl geometry in 3D, self-dual conformal geometry in 4D on every solution).
\item Integrability {\it on equation}, meaning  integrability (in twistor-theoretic sense) of the canonical $GL(2, \mathbb{R})$ structure induced on a fourfold $X\subset {\bf Gr}(3, 5)$.
\end{itemize}

We demonstrate that for $(d, n)=(3, 5)$ these seemingly different approaches lead to one and the same class of integrable systems $\Sigma(X)$ and
prove that the parameter space of such systems is 30-dimensional. Factored by the natural action of the equivalence group ${\bf SL}(5)$, this gives 6-dimensional moduli space of integrable systems $\Sigma(X)$. We give a complete description of linearisable systems (the corresponding fourfold $X$ is a linear section of ${\bf Gr}(3, 5)$) and linearly degenerate systems (the corresponding fourfold $X$ is the image of a quadratic map $\mathbb{P}^4\dashrightarrow {\bf Gr}(3, 5)$). The fourfolds corresponding to `generic' integrable systems are not algebraic, and can be parametrised by generalised hypergeometric functions.


\bigskip


\bigskip

\end{abstract}

\newpage

\tableofcontents

\newpage

\section{Introduction}

\subsection{Formulation of the problem}

In the most general setting the problem that we address in this paper can be described as follows. Let ${\bf Gr}(d, n)$ be the Grassmannian of $d$-dimensional linear subspaces of an $n$-dimensional vector space $V^n$. A submanifold $X\subset {\bf Gr}(d, n)$ gives rise to the differential system $\Sigma(X)$ that governs $d$-dimensional submanifolds of $V^n$ whose Gaussian image is contained in $X$ (we recall that the Gaussian image of a  submanifold  is the collection of its tangent spaces translated to the origin). In this sense we are in the context of Grassmann geometries as discussed in \cite{Harvey}.
Since  $d$-dimensional submanifolds of $V^n$ are (locally) parametrised by  $n-d$ functions of $d$ variables,  we will assume that the codimension of $X$ in  ${\bf Gr}(d, n)$ also equals $n-d$: in this case $\Sigma(X)$ will be a determined system of $n-d$ first-order PDEs for $n-d$ unknown functions of $d$ independent variables. Systems of  the form $\Sigma(X)$ appear in a wide range of applications in differential geometry (in particular, calibrated geometries \cite{Harvey}),  
general relativity (Einstein-Weyl structures \cite{D}, heavenly-type equations \cite{Plebanski}), and the theory of integrable systems (B\"acklund transformations in higher dimensions, dispersionless Lax pairs \cite{Zakharov}, dispersionless limits of multi-dimensional  soliton equations). Particularly interesting
examples are provided by multidimensional ($d\geq 3$) {\it integrable} systems $\Sigma(X)$: such systems possess an infinity of multi-phase solutions that can be considered as analogues of multi-gap solutions of integrable soliton equations. It was observed in \cite{Fer1} that the requirement of the existence of such solutions constitutes an efficient  criterion (known as the method of hydrodynamic reductions: see Sect. \ref{sec:method} for a brief description of the method),  that allows one to  classify integrable systems of the form $\Sigma(X)$. In this paper we concentrate on the following key questions:

\begin{itemize}

\item {What is the dimension of the moduli space of integrable systems $\Sigma(X)$? What are the most interesting examples?}

\item {How does integrability of a system $\Sigma(X)$  translate into differential geometry of the corresponding submanifold $X\subset {\bf Gr}(d, n)$?}

\end{itemize}

In the case $n=d+1$ the system $\Sigma(X)$ consists of a single first-order PDE for a scalar function of $d$ independent variables, which can be solved by the method of characteristics via reduction to ODEs.  Thus,  we  can assume $ n\geq d+2$. The case $d=1$ also corresponds to an ODE system. 
The case $d=2$, $n$ arbitrary,   leads to another familiar class of systems: introducing in $V^n$ coordinates $t, x, u^{i}, \ i=1, \dots, n-2$, and parametrising two-dimensional submanifolds of $V^n$ in the form $u^i=u^i(t, x)$, we can represent the corresponding  system $\Sigma(X)$ as
$$
F^1(u^{i}_t,  u^{i}_x)=0,  ~~ \dots, ~~ F^{n-2}(u^{i}_t,  u^{i}_x)=0.
$$
Solving these equations for $u^{i}_t$ in the form
$u^{i}_t=f^{i}({\bf u}_x),$ here ${\bf u}=(u^1, \dots, u^{n-2})$,
differentiating  by $x$ and setting ${\bf u}_x={\bf v}$, we obtain a  system of conservation laws,
$$
v^{i}_t=\partial_{x} f^{i}({\bf v}).
$$
Differential-geometric and integrability aspects of such systems were extensively studied in \cite{Tsar, Sevennec, Dubrovin, Agafonov}, see also references therein. Thus, in what follows we  assume $n\geq d+2$ and $d\geq 3$.

\medskip

To be more specific we will concentrate on the  particularly interesting case of this problem  where our results are fairly complete, the case of fourfolds $X\subset {\bf Gr}(3, 5)$ ($n=5, \ d=3$). Introducing in $V^5$ coordinates $x^1, x^2, x^3, u, v,$ one  can parametrise three-dimensional submanifolds of $V^5$ in the form $u=u(x^1, x^2, x^3), \ v=v(x^1, x^2, x^3)$. Their tangent spaces are given by $du=u_idx^i, \ dv=v_idx^i$ where  $u_i, \ v_i$ can be viewed as local coordinates on ${\bf Gr}(3, 5)$.
The corresponding system $\Sigma(X)$ reduces to a pair of first-order PDEs for $u$ and $v$,
\begin{equation}
F(u_{1}, u_{2}, u_{3}, v_{1}, v_{2}, v_{3})=0, ~~~ G(u_{1}, u_{2}, u_{3}, v_{1}, v_{2}, v_{3})=0,
\label{main}
\end{equation}
$u_i=\partial u/\partial x^i, \  v_i=\partial v/\partial x^i$. Equations (\ref{main})  specify a fourfold  $X\subset {\bf Gr}(3, 5)$. The class of systems (\ref{main}) is invariant under the equivalence group ${\bf SL}(5)$ that acts by linear transformations on the combined set of variables $ x^1, x^2, x^3, u, v$. Since this action preserves the integrability, all our classification results will be formulated modulo ${\bf SL}(5)$-equivalence. Necessary details on the equivalence group are provided in Sect. \ref{sec:equiv}.

There exists a whole variety of important examples that fall into  class (\ref{main}). One of them appears in the context of the  dispersionless Kadomtsev-Petviashvili (dKP) equation, $u_{xt}-u_xu_{xx}-u_{yy}=0$,  one of the most well-studied  dispersionless integrable PDEs arising in nonlinear acoustics \cite{KZ} and the theory of Einstein-Weyl structures \cite{D}. Its dispersionless Lax pair \cite{Zakharov}
consists of two first-order relations of type (\ref{main}),
\begin{equation}
v_y-\frac{1}{2}v_x^2-u_x=0, ~~~ v_t-\frac{1}{3}v_x^3-v_xu_x-u_y=0,
\label{dKP}
\end{equation}
here  $(x^1, x^2, x^3)=(x, y, t)$. The dKP equation results from (\ref{dKP}) on  elimination of $v$, that is, via the compatibility condition $v_{yt}=v_{ty}$. Similarly, the elimination of $u$ leads to the modified dKP (mdKP) equation, $v_{xt}-(v_y-\frac{1}{2}v_x^2)v_{xx}-v_{yy}=0.$ Thus,  relations (\ref{dKP}) provide   B\"acklund-type transformation connecting  dKP and mdKP equations.  We refer to Sect. \ref{sec:ex} for further examples and classification results.

\subsection{Non-degeneracy condition}
\label{sec:nondeg}

In what follows we assume that system (\ref{main}) is {\it non-degenerate} in
the following sense. Let us consider
the corresponding linearised system,
$$
\left(\begin{array}{cc}
F_{u_1} & F_{v_1}\\
G_{u_1} & G_{v_1}
\end{array}\right)
\left(
\begin{array}{c}
\cal U\\
\cal V
\end{array}\right)_{x^1}+
\left(\begin{array}{cc}
F_{u_2} & F_{v_2}\\
G_{u_2} & G_{v_2}
\end{array}\right)
\left(
\begin{array}{c}
\cal U\\
\cal V
\end{array}\right)_{x^2}+
\left(\begin{array}{cc}
F_{u_3} & F_{v_3}\\
G_{u_3} & G_{v_3}
\end{array}\right)
\left(
\begin{array}{c}
\cal U\\
\cal V
\end{array}\right)_{x^3}=0,
$$
obtained by setting $u\to u+\epsilon\, {\cal U}, \ v\to v+\epsilon {\cal V}$, expanding $F$ and $G$ in Taylor series  and keeping  terms of the
order $\epsilon$. Non-degeneracy means that the dispersion relation (characteristic variety),
$$
\det \left[
\lambda^1 \left(\begin{array}{cc}
F_{u_1} & F_{v_1}\\
G_{u_1} & G_{v_1}
\end{array}\right)
+
\lambda^2 \left(\begin{array}{cc}
F_{u_2} & F_{v_2}\\
G_{u_2} & G_{v_2}
\end{array}\right)
+
\lambda^3 \left(\begin{array}{cc}
F_{u_3} & F_{v_3}\\
G_{u_3} & G_{v_3}
\end{array}\right)
\right]=0,
$$
 defines an
irreducible conic in $\mathbb{P}^2$ with homogeneous coordinates  $(\lambda^1: \lambda^2: \lambda^3)$. Explicitly, the dispersion relation can be represented in the form
$(\lambda^1, \lambda^2, \lambda^3)g^{\sharp}(\lambda^1, \lambda^2, \lambda^3)^t=0$ where $g^{\sharp}$ is the $3\times 3$ symmetric matrix,
$$
g^{\sharp}=g^{ij}=\frac{1}{2}(F_{u^i}G_{v^j}+F_{u^j}G_{v^i}-F_{v^i}G_{u^j}-F_{v^j}G_{u^i}).
$$
It gives rise to the conformal structure
$g=g_{ij}dx^idx^j$ (here $g_{ij}$ is the inverse of $g^{ij}$).   Note that  non-degeneracy  is equivalent to $\det g \ne 0$:  this is the case for all known systems of physical/geometric relevance.  
It turns out that the signature of $g$ is always  Lorentzian, and thus our PDE system is hyperbolic.

\com{
For four-dimensional  systems,
\begin{equation}
 F(u_i,  v_i)=0, ~~~ G(u_i, v_i)=0,
\label{d}
\end{equation}
where $u$ and $v$ are functions of  $x^1, \dots, x^4$,
 non-degeneracy means that the dispersion relation,
$$
\det \left[\sum_{i=1}^4
\lambda^i \left(\begin{array}{cc}
F_{u_i} & F_{v_i}\\
G_{u_i} & G_{v_i}
\end{array}\right)
\right]=0,
$$
defines an
irreducible quadric of rank four. This quadric gives rise to a non-degenerate  conformal structure  in 4D. }

Geometric aspects  of  conformal structures defined by the characteristic variety
will play a  key role  in our characterisation of integrable systems: we will see that solutions to  integrable equations carry `integrable'  background geometry. In 3D, this is the Einstein-Weyl geometry.

\subsection{Einstein-Weyl geometry in 3D}

Recall that an Einstein-Weyl structure  consists of a symmetric connection $\mathbb{D}$ and a conformal structure $g$  such that:

\noindent (a) connection $\mathbb {D}$ preserves the conformal class of $g$: $\mathbb{D}[g]=0$,

\noindent (b)  trace-free part of the symmetrized Ricci tensor of $\mathbb {D}$ vanishes.

In coordinates, this gives
\begin{equation}
\mathbb{D}_kg_{ij}=\omega_k g_{ij}, ~~~ R_{(ij)}=\Lambda g_{ij},
\label{EW}
\end{equation}
where $\omega=\omega_kdx^k$ is a covector, $R_{(ij)}$ is the symmetrized Ricci tensor of $\mathbb{D}$,  and $\Lambda$ is some function \cite{Cartan}.  Note that it is sufficient  to specify $g$ and $\omega$ only, then the first set of equations uniquely determines $\mathbb{D}$. The integrability of Einstein-Weyl equations  (\ref{EW}) by   twistor-theoretic methods was established by Hitchin  \cite{Hitchin}. It was shown in \cite{DFK} that generic Einstein-Weyl structures are governed by the Manakov-Santini system introduced  in \cite{Man-San} as a two-component integrable generalisation of the dKP equation.
We will see that solutions to integrable systems
(\ref{main}) carry Einstein-Weyl geometry: the conformal structure $g$ defined by the characteristic variety  must be Einstein-Weyl on every solution,  furthermore,  the covector $\omega$ can be expressed in terms of  $g$ by the universal  formula
\begin{equation}
\omega_k=2g_{kj}\mathcal{D}_{x^s}(g^{js})+\mathcal{D}_{x^k}(\ln\det g_{ij}),
\label{omega}
\end{equation}
where $\mathcal{D}_{x^k}$ denotes  total derivative with  respect to the independent variable $x^k$ (note that $g$ depends on first-order jets of the solution $u, v$). Formula  (\ref{omega}) appeared  in \cite{FerKrug} in geometric approach to the dispersionless integrability in 3D. It is  invariant
under the gauge transformation $g\to \lambda g, \ \omega \to \omega + d\ln\lambda$,
the property characteristic of  Einstein-Weyl geometry.
 According to the result of Cartan \cite{Cartan}, the Einstein-Weyl property of a triple  $(\mathbb{D}, \ g, \ \omega)$ is equivalent to the existence of a two-parameter family of surfaces that are null with respect to the conformal structure $g$ (that is, tangential to the null cones of $g$), and totally geodesic in the Weyl connection  $\mathbb{D}$. In the context of  integrable systems (\ref{main}), such surfaces are provided by the corresponding dispersionless Lax pairs: these consist  of  $\lambda$-dependent vector fields $X, Y$ that are required to commute modulo (\ref{main}),  identically in the `spectral parameter' $\lambda$. For systems (\ref{main}),  the existence of such Lax pairs is equivalent to the integrability by the method of hydrodynamic reductions,  see Sect. \ref{sec:Lax}.  Taking integral surfaces of the distribution spanned by $X, Y$ in the extended four-space with coordinates $x, y, t, \lambda$, and projecting them down to the space of independent variables $x, y, t$, we obtain the required two-parameter family of null totally geodesic surfaces.
 Let us mention that  relations of dispersionless integrable systems in 3D to
Einstein-Weyl geometry have been discussed previously in \cite{Ward, Calderbank,   Dun7, FerKrug}, see also references therein.

\noindent {\bf Example 1.} Let us consider the system
$$
u_{t}-\frac{1}{2}u_x^2-v_y=0, ~~~ v_x-u_y=0,
$$
which reduces to the dKP equation, $u_{xt}-u_xu_{xx}-u_{yy}=0$, on elimination on $v$. Its characteristic variety  defines the conformal structure
$g=4dxdt-dy^2+4u_xdt^2$. Introducing the covector $\omega=-4u_{xx}dt$ by formula (\ref{omega}), one can verify that the pair $g, \omega$ satisfies the Einstein-Weyl equations if  $u$ solves the dKP equation.  This Einstein-Weyl structure was obtained previously in \cite{D}. The corresponding Lax pair has the form
$$
X=\partial_y-\lambda \partial_x +u_{xx}\partial_{\lambda}, ~~~
Y=\partial_t-(\lambda^2+u_x) \partial_x +(u_{xx}\lambda +u_{xy})\partial_{\lambda};
$$
one can verify that these vector fields commute modulo the above system. Projecting integral surfaces of the distribution spanned by $X, Y$ from the extended space of variables $x, y, t, \lambda$ to the space of independent variables $x, y, t$, one obtains a two-parameter family of null totally geodesic surfaces of the corresponding Einstein-Weyl structure.

\subsection{${\bf GL}(2, \bbbr)$ structures}
\label{sec:Bryant}

The  tangent bundle to the Grassmannian  ${\bf Gr}(3, 5)$ carries  canonical generalised conformal structure defined by the family of Segre cones  $du_idv_j-du_jdv_i=0$. Thus, each projectivised tangent space ${\mathbb{P}}{\rm T} {\bf Gr}(3, 5)$ contains a Segre variety,  an algebraic threefold of degree three. Given a non-degenerate fourfold $X\subset {\bf Gr}(3, 5)$
(a fourfold is said to be non-degenerate if it gives rise to a non-degenerate system (\ref{main})), the intersection of its tangent space ${\rm T}X$ with the Segre cone is a two-dimensional rational cone of degree three; its projectivisation is a rational normal curve of degree three (twisted cubic). Thus,  ${\mathbb{P}}{\rm T}X$ is supplied with a field of twisted cubics. This is known as a ${\bf GL}(2, \bbbr)$ structure on $X$. It was demonstrated by Bryant \cite{Bryant} that every  four-dimensional  ${\bf GL}(2, \bbbr)$ structure defines on $X$ a canonical affine connection which preserves the ${\bf GL}(2, \bbbr)$ structure, and  whose torsion lies in  $8$-dimensional irreducible representation of  ${\bf GL}(2, \bbbr)$. We will call it the Bryant connection, see Sect. \ref{sec:curvtor} for computational formulae. Various important properties of system (\ref{main}) have natural interpretation in terms of this connection:
\begin{itemize}


\item System (\ref{main}) is linearly degenerate and integrable if and only if the  Bryant connection is symmetric and flat (Proposition 5 of Sect. \ref{sec:lindeg}).

\item System (\ref{main}) is integrable if and only if the curvature $R$ and the covariant derivative $\nabla T$ of the torsion $T$ of the   Bryant connection  are certain invariant quadratic expressions in $T$,
$$
R=f(T^2), ~~~ \nabla T=g(T^2),
$$
see Theorem \ref{decomp} of Sect. \ref{sec:curvtor} for precise statements. These expressions are analogous to the ones obtained by Smith \cite{Smith} in the context of 5-dimensional ${\bf GL}(2, \bbbr)$ structures associated with integrable equations of the dispersioless Hirota type \cite{Hirota}. Our expressions provide a compact invariant formulation of the integrability conditions.

\end{itemize}


\subsection{The method of hydrodynamic reductions}
\label{sec:method}

In the most general set-up,  the method of hydrodynamic reductions \cite{Fer1} applies to quasilinear systems of the form
\begin{equation}
A ({\bf u}){\bf u}_x+B({\bf u}){\bf u}_y+C({\bf u}){\bf u}_t=0,
\label{quasi1}
\end{equation}
where ${\bf u}=(u^1, ..., u^m)^t$ is an $m$-component column vector of the dependent variables, and $A, B, C$ are $l\times m$ matrices where $l$, the number of equations, is allowed to exceed the number of  unknowns, $m$. Note that system (\ref{main}) can be brought to  form  (\ref{quasi1}) by choosing $u_i, v_i$ as the new dependent variables and writing out all possible consistency conditions among them - see Sect. \ref{sec:hydro} for  details. The method of hydrodynamic reductions consists of seeking multiphase solutions in the form
$$
{\bf u}={\bf u}(R^1, ..., R^N),
$$
where the phases $R^i(x, y, t)$, whose number $N$ is allowed to be arbitrary, are required to satisfy a pair of consistent $(1+1)$-dimensional
systems,
\begin{equation}
R^i_t=\lambda^i(R) R^i_x, ~~~ R^i_y=\mu^i(R) R^i_x,
\label{R}
\end{equation}
known as systems of hydrodynamic type. The corresponding characteristic speeds $\lambda^i$ and $\mu^i$ are required to satisfy the commutativity conditions \cite{Tsar},
\begin{equation}
\frac{\partial_j\lambda
^i}{\lambda^j-\lambda^i}=\frac{\partial_j\mu^i}{\mu^j-\mu^i},
\label{comm}
\end{equation}
here $i\ne j, \  \partial_j=\partial_{ R^j}$.
Solutions of this type originate from gas dynamics, and are  known as nonlinear interactions of planar simple waves. They  can be interpreted as natural dispersionless analogues of finite-gap solutions of $2+1$ dimensional soliton equations. Equations (\ref{R}) are said to define an $N$-component  hydrodynamic reduction of the original system (\ref{quasi1}). System (\ref{quasi1}) is said to be {\it integrable} if, for every $N$, it possesses infinitely many $N$-component hydrodynamic reductions parametrised by $N$ arbitrary functions of one variable. This requirement imposes strong constraints (integrability conditions) on the matrix elements of  $A, B$ and $C$. We refer to \cite{Fer1} for further details and references.

\subsection{Summary of the main results}

In Sect. \ref{sec:hydro}
 we apply the  method of hydrodynamic reductions to  equations of type (\ref{main}). This results in an overdetermined involutive system of PDEs for the functions $F$ and $G$ (the integrability conditions). The analysis of this system  leads to our first  result (see Theorem \ref{thm1} of section \ref{sec:hydro} and Proposition 2 of section \ref{sec:equiv}):
 
 \medskip
 
{\it 
The moduli space of non-degenerate integrable systems  (\ref{main}) is $6$-dimensional. Furthermore, the following conditions are equivalent:

\noindent {(a)} System (\ref{main}) is integrable by the method of hydrodynamic reductions.

\noindent  {(b)} Conformal structure $g$  defined by the characteristic variety of system (\ref{main}) is Einstein-Weyl on every solution (with  covector $\omega$  given by formula (\ref{omega})).

}

\medskip






Thus, Einstein-Weyl equations (\ref{EW}) provide an efficient integrability test. Using the integrability conditions we prove that every integrable system (\ref{main}) possesses a dispersionless Lax pair (Theorem \ref{thm2} of Sect. \ref{sec:Lax}).   Furthermore, we describe a construction that links  Lax pairs  to null totally geodesic surfaces of the corresponding Einstein-Weyl structure (Sect. \ref{sec:null}).

In Sect. \ref{sec:equiv} we demonstrate that the class of equations (\ref{main}) is invariant under the equivalence group  ${\bf SL}(5)$  that acts by linear transformations on the combined set of variables $ x^1, x^2, x^3, u, v$. All our results are formulated modulo this equivalence.

In Sect. \ref{sec:ex} we  describe a variety of   integrable examples of type (\ref{main}) expressible in elementary functions, modular forms and theta functions. In particular, we demonstrate that system (\ref{main}) is linearisable by a transformation from the equivalence group if and only if the conformal structure defined by its characteristic variety is  flat on every solution (Sect. \ref{sec:lin}). Next, we show that linearly degenerate integrable systems are characterised by the property that the associated Bryant connection is symmetric and flat. Fourfolds $X$ corresponding to such systems are images of quadratic maps  $\mathbb{P}^{4}\dashrightarrow {\bf Gr}(3, 5)$ given by the classical construction of Chasles. This provides a complete list of normal forms of linearly degenerate integrable systems (Sect. \ref{sec:lindeg}).
Based on the  construction of Odesskii and Sokolov \cite{Odesskii2, Odesskii3} we point out that   {\it generic} integrable systems (\ref{main}) can be parametrised by generalised hypergeometric functions (Sect. \ref{sec:hyper}). This  demonstrates that  fourfolds $X\subset {\bf Gr}(3, 5)$ corresponding to integrable systems (\ref{main}) may have a fairly intricate analytic structure: in particular, they do not need to be  algebraic. Among the simplest nontrivial  examples of this kind one can mention the system
$$
u_t=v_x, ~~~ v_{t}=\frac{u_{y}}{v_{x}}+\frac{1}{6}\eta (u_{x})v_{x}^{2},
$$
which appeared in  \cite{MaksEgor} in the classification of integrable hydrodynamic chains. In this case the integrability conditions result in the Chazy equation for $\eta$,  
$
\eta^{\prime \prime \prime }+2\eta \eta ^{\prime \prime }-3{\eta^{\prime} }^2=0
$
(here prime denotes differentiation by $u_x$), whose generic solution is  expressible in terms of the  Eisenstein series $E_2$.

 In Sect. \ref{sec:geom} we demonstrate that every fourfold $X\subset {\bf Gr}(3, 5)$ inherits a ${\bf GL}(2, \bbbr)$ structure, namely, a field of twisted cubics specified in the projectivised tangent bundle $\mathbb{P}{\rm T}X$ (Sect. \ref{sec:gl2}).
 The integrability can be reformulated geometrically as the requirement that the associated ${\bf GL}(2, \bbbr)$ structure possesses infinitely many three-dimensional submanifolds carrying a holonomic characteristic net (Sect. \ref{sec:secant}). We  show that every fourfold $X$ corresponding to integrable system (\ref{main}) carries  canonical conformally symplectic structure  which is parallel in the Bryant connection (Sect. \ref{sec:confsymp}).
Finally, we reformulate the integrability conditions of system (\ref{main}) in terms of the curvature and torsion of the associated  connection (Theorem \ref{decomp} of Sect. \ref{sec:curvtor}).
   We emphasize that ${\bf GL}(2, \bbbr)$ structures arising in this context are quite different from those discussed in \cite{Bryant}, in particular, they  generally have nonzero torsion.

Throughout the paper, most of our considerations are local, although some results related to algebraic fourfolds $X\subset  {\bf Gr}(3, 5)$ have global nature. 
We are primarily interested in real smooth structures and hyperbolic systems, however, the theory applies to the complex-analytic case as well. 
In the calculations we use computer algebra systems Mathematica and Maple (these only utilise symbolic polynomial algebra over $\mathbb{Q}$, so the results are rigorous). The programs are available from 	arXiv:1503.02274v1 [math.DG].

\section{Integrability conditions and Lax pairs }
\label{sec:IL}

\subsection{Derivation of  integrability conditions}
\label{sec:hydro}

In this section we establish our first main result:

\begin{theorem}\label{thm1}
{\it  
The parameter space of non-degenerate integrable systems  (\ref{main}) is $30$-dimensional.
 Furthermore, the following conditions are equivalent:

\noindent {(a)} System (\ref{main}) is integrable by the method of hydrodynamic reductions.

\noindent  {(b)} Conformal structure $g$  defined by the characteristic variety of system (\ref{main}) is Einstein-Weyl on every solution (with  covector $\omega$  given by formula (\ref{omega}))}.


\end{theorem}

\centerline {\bf Proof:}

\medskip
  \noindent

 Applied to system (\ref{main}),  the method of hydrodynamic reductions leads to a set of differential constraints (integrability conditions) for the functions $F$ and $G$. The 30-dimensionality of the parameter space follows from the involutivity of these  conditions.  Once the integrability conditions are derived, the equivalence of (a) and (b) can  be shown  by a direct calculation. Alternatively, one can utilise the fact that every integrable system (\ref{main}) possesses a dispersionless Lax pair (Theorem \ref{thm2} of Sect. \ref{sec:Lax}), which provides null totally geodesic surfaces of the corresponding  Einstein-Weyl structure $(\mathbb{D}, g, \omega)$. This implies the Einstein-Weyl property due to the result of Cartan \cite{Cartan}.


To derive the integrability conditions we proceed  as follows.
Let us  first rewrite  system  (\ref{main}) in  evolutionary form,
\begin{equation}
u_{t}=f(u_{x}, u_{y}, v_{x}, v_{y}), ~~~ v_{t}=g(u_{x}, u_{y}, v_{x}, v_{y}),
\label{evol}
\end{equation}
here we use the notation $(x^1, x^2, x^3)=(x, y, t)$. Our strategy is to derive a set of constraints for the right hand sides $f$ and $ g$ that are necessary and sufficient for  integrability. We begin by  transforming our system into  first-order quasilinear form. This can be done by taking first-order partial derivatives of $u$ and $v$ as the new dependent variables and writing out all possible consistency conditions among them. Thus, we introduce the notation
$$
u_{x}=a,\ u_{y}=b, \ v_{x}=p,\ v_{y}=q, \ u_{t}= f(a,b,  p, q), \ v_t=g(a, b, p, q).
$$
This provides an equivalent  quasilinear representation,
\begin{equation}
\begin{array}{c}
a_y=b_x, ~~ a_t=f(a, b, p, q)_x, ~~ b_t=f(a, b, p, q)_y, \\
\ \\
p_y=q_x,  ~~ p_t=g(a,b,  p, q)_x, ~~ q_t=g(a,b,  p, q)_y.
\end{array}
\label{quasi}
\end{equation}
Note that   quasilinear representation (\ref{quasi}) is of type (\ref{quasi1}) with  $m=4,  \ l=6, \ {\bf u}=(a, b, p, q)$.
Looking for multi-phase solutions in the form
$$
a=a(R^{1},...,R^{N}),\ b=b(R^{1},...,R^{N}), \ p=p(R^{1},...,R^{N}), \ q=q(R^{1},...,R^{N}),
$$
 where the phases  $R^i$ satisfy   equations (\ref{R}),
and substituting this ansatz into (\ref{quasi}), we obtain the relations
\begin{equation}
\partial_ib=\mu^{i}\partial_ia, ~~  \partial_iq=\mu^{i}\partial_ip,
\label{E1}
\end{equation}
here $ \partial_i=\partial_{ R^i}$, as well as
\begin{equation}
\begin{array}{c}
(\lambda^i-f_a-\mu^if_b)\partial_ia=(f_p+\mu^if_q)\partial_ip, ~~~
(\lambda^i-g_p-\mu^ig_q)\partial_ip=(g_a+\mu^ig_b)\partial_ia.
\end{array}
\label{E2}
\end{equation}
The last two equations imply  the dispersion relation connecting $\lambda^i$ and $\mu^i$,
\begin{equation}
(\lambda^i-f_a-\mu^if_b)(\lambda^i-g_p-\mu^ig_q)=(f_p+\mu^if_q)(g_a+\mu^ig_b).
\label{conic}
\end{equation}
In what follows we assume that the dispersion relation defines an irreducible conic in the $(\lambda, \mu)$-plane: this is equivalent to    the non-degeneracy of  system (\ref{main}) as discussed in Sect. \ref{sec:nondeg}.   Setting in (\ref{E2}) $\partial_i a=\varphi^i\partial_ip$
and solving the resulting equations for $\lambda^i$ and $\mu^i$
we can parametrise  dispersion relation (\ref{conic}) in the form
$$
\mu^i=-\frac{f_p+(f_a-g_p)\varphi^i-g_a{\varphi^i}^2}{f_q+(f_b-g_q)\varphi^i-g_b{\varphi^i}^2}, ~~~
\lambda^i=\frac{(f_q+f_b\varphi^i)(g_p+g_a \varphi^i)-(f_p+f_a\varphi^i)(g_q+g_b \varphi^i)}{f_q+(f_b-g_q)\varphi^i-g_b{\varphi^i}^2}.
$$
Substituting these parametric expressions into the commutativity conditions (\ref{comm}), and using the relations
\begin{equation}
\partial_ia=\varphi^{i}\partial_ip, ~~~ \partial_ib=\mu^{i}\varphi^i\partial_ip, ~~  \partial_iq=\mu^{i}\partial_ip,
\label{abq}
\end{equation}
we obtain   $\partial_j \varphi^i$ in the form
$
\partial_j \varphi^i=(\dots) \partial_jp, \ i\ne j,
$
where dots denote  rational expressions in $\varphi^i, \ \varphi^j$ whose coefficients depend on second-order partial derivatives of $f$ and $g$ (we omit these expressions due to their complexity). Calculating  consistency conditions for relations (\ref{abq}), $\partial_i\partial_ja=\partial_j\partial_ia, \ \partial_i\partial_jb=\partial_j\partial_ib, \ \partial_i\partial_jq=\partial_j\partial_iq,$ we obtain (one and the same!) expression for $\partial_i\partial_jp$ in the form $\partial_i\partial_jp=(\dots)\partial_ip\partial_jp, \ i\ne j,$ where, again, dots denote terms rational in $\varphi^i$ and $\varphi^j$.
Ultimately,   $N$-phase solutions are governed by the relations
\begin{equation}
\partial_j \varphi^i=(\dots) \partial_jp, ~~~ \partial_i\partial_jp=(\dots)\partial_ip\partial_jp, ~~~ i\ne j.
\label{int}
\end{equation}
To ensure their solvability we need to impose the compatibility conditions
\begin{equation}
\partial_k\partial_j \varphi^i=\partial_j\partial_k \varphi^i, ~~~
\partial_k\partial_j\partial_ip=\partial_j\partial_k\partial_ip,
\label{comp}
\end{equation}
that are required to hold for every triple of indices $i\ne j \ne k$. Direct calculation based on (\ref{abq}) and (\ref{int}) results in
$$
\partial_k\partial_j \varphi^i-\partial_j\partial_k \varphi^i=(\dots)\partial_jp\partial_kp, ~~~
\partial_k\partial_j\partial_ip-\partial_j\partial_k\partial_ip=(\dots)\partial_ip\partial_jp\partial_kp,
$$
where dots denote  rational expressions in $\varphi^i, \varphi^j, \varphi^k$ whose coefficients depend on partial derivatives of $f$ and $g$ up to  order three. To ensure  solvability of equations (\ref{int}) we set all these coefficients equal to zero. This gives a system of differential constraints for $f$ and $g$ that is  linear in the third-order derivatives thereof.
Moreover, the number of linearly independent compatibility conditions equals  the total number  of third-order partial derivatives of $f$ and $g$, namely $40$.  Solving for third-order derivatives of $f$ and $g$ we obtain the required integrabililty conditions that can be represented in symbolic form,
\begin{equation}
d^3f=R(df, dg, d^2f, d^2g), ~~~ d^3g=S(df, dg, d^2f, d^2g),
\label{*}
\end{equation}
$40$ equations altogether. Here $R$ and $S$ depend  rationally on the first- and second-order partial derivatives of $f$ and $g$.  Straightforward  calculation shows that overdetermined system (\ref{*}) is in involution. Thus, the parameter space of integrable systems  (\ref{main}) is $30$-dimensional: one can arbitrarily specify the values of $f, g$, $df, dg$ and $d^2f, d^2g$ at any fixed point. This amounts to $2\times 1+2\times 4+2\times 10=30$ arbitrary constants.

\medskip

Finally, the equivalence of (a) and (b) can be established by a direct calculation: in coordinates $(x, y, t)$, the characteristic variety of system (\ref{evol}) defines contravariant metric (metric with upper indices) with the matrix
$$
g^{ij}=\left(
\begin{array}{ccc}
f_ag_p-f_pg_a & \frac{1}{2}(f_ag_q-f_qg_a+f_bg_p-f_pg_b) & -\frac{1}{2}(f_a+g_p)\\
\ \\
 \frac{1}{2}(f_ag_q-f_qg_a+f_bg_p-f_pg_b) & f_bg_q-f_qg_b& -\frac{1}{2}(f_b+g_q) \\
\ \\
 -\frac{1}{2}(f_a+g_p) &  -\frac{1}{2}(f_b+g_q) & 1
 \end{array}
 \right).
$$
Introducing $\omega$ by formula (\ref{omega}) and calculating  Einstein-Weyl conditions (\ref{EW}) we obtain expressions depending on third-order partial derivatives of $u$ and $v$. These have to vanish identically modulo (\ref{evol}), in other words,  on every solution. Using (\ref{evol}) and its differential consequences to eliminate all partial derivatives involving differentiation by $t$,  we obtain differential expressions that are polynomial in the remaining second- and third-order partial derivatives of $u$ and $v$. It can be verified directly that the vanishing of all coefficients of these polynomials is equivalent to  integrability conditions (\ref{*}). More precisely, the Einstein-Weyl conditions contain terms of two types: linear in third-order derivatives, and quadratic in second-order derivatives  of $u$ and $v$. Coefficients at third-order derivatives of $u$ and $v$ vanish identically due to the choice of $\omega$, while coefficients at quadratic terms (105 coefficients altogether) give all of the 40 integrability conditions (\ref{*}).
Details of calculations described above, including the integrability conditions,  are available from arXiv:1503.02274v1 [math.DG]. 
 This finishes the  proof of Theorem \ref{thm1}.




\medskip

\noindent {\bf Remark.} Note that each compatibility condition (\ref{comp})  involves three distinct  indices only. This leads to the following important observations.

\noindent {\bf (i)} Every system of type (\ref{main}) has infinitely many  two-component  reductions  parametrised by two arbitrary functions of one variable. Indeed, these reductions are governed by  equations (\ref{int})  where $i, j=1, 2$. Since in this case equations (\ref{int}) are automatically consistent,  their general solution depends, modulo reparametrisations $R^i\to \varphi^i(R^i)$, on two arbitrary functions of one variable.
Therefore, the existence of two-component reductions is a common phenomenon that is not related to the  integrability.

\noindent {\bf (ii)}  On the contrary, the existence of three-component reductions implies the existence of $N$-component reductions for arbitrary $N$. Thus, one can define integrability as the existence of infinitely many three-component reductions parametrised by three arbitrary functions of one variable. This property is reminiscent of the well-known three-soliton condition in the theory of $(1+1)$-dimensional integrable soliton equations.

 \medskip


\subsection{Dispersionless Lax pairs}
\label{sec:Lax}

We say that system (\ref{main}) possesses a dispersionless Lax pair if there exists an auxiliary function $S$ satisfying two   Hamilton-Jacobi type equations,
\begin{equation}
S_2=P(S_1, u_i, v_i), ~~~ S_3=Q(S_1, u_i, v_i),
\label{Lax}
\end{equation}
such that the compatibility condition, $S_{23}=S_{32}$,  holds identically modulo (\ref{main}); here $S_i=\partial S/\partial {x^i}$. Note that the dependence of $P$ and $Q$ on $S_1$ is usually nonlinear.   Lax pairs of this type appeared in the context of dispersionless integrability in \cite{Zakharov}; in many cases, including the one studied in this paper, the existence of a dispersionless Lax pair is equivalent to the integrability by the method of hydrodynamic reductions. For instance, system ({\ref{dKP}) possesses the Lax pair
$$
S_y=S_x^2+v_xS_x, ~~~ S_t=\frac{4}{3}S_x^3+2v_xS_x^2+(u_x+v_x^2)S_x;
$$
the compatibility condition $S_{yt}=S_{ty}$ is satisfied identically modulo (\ref{dKP}). The main result of this section is the following:

 \begin{theorem}\label{thm2} Integrability of  system (\ref{main})  is equivalent to the existence  of a `generic' dispersionless Lax pair.
 \end{theorem}

\medskip

\noindent The definition of the term `generic' is given after formula (\ref{PQ}).
\medskip

\centerline  {\bf Proof:}

\medskip
 \noindent Let us  again consider system (\ref{main}) in evolutionary form,
$$
u_{t}=f(u_{x}, u_{y}, v_{x}, v_{y}), ~~~ v_{t}=g(u_{x}, u_{y}, v_{x}, v_{y}).
$$
Setting $u_{x}=a,\ u_{y}=b, \ v_{x}=p,\ v_{y}=q$ we will look for a Lax pair in the form
\begin{equation}
S_y=P(S_x, a, b, p, q), ~~~ S_t=Q(S_x, a, b, p, q).
\label{Lax1}
\end{equation}
Requiring that  the consistency condition $S_{yt}=S_{ty}$ holds identically modulo (\ref{quasi}), we obtain $6$ first-order relations for $P$ and $Q$:
\begin{equation}
\begin{array}{c}
f_aP_a+g_aP_p+P_{\lambda}Q_a-Q_{\lambda}P_a=0, ~~~ f_pP_a+g_pP_p+P_{\lambda}Q_p-Q_{\lambda}P_p=0,\\
\ \\
Q_a=f_bP_a+f_aP_b+g_bP_p+g_aP_q+P_{\lambda}Q_b-Q_{\lambda}P_b,\\
\ \\
Q_p=f_qP_a+f_pP_b+g_qP_p+g_pP_q+P_{\lambda}Q_q-Q_{\lambda}P_q,\\
\ \\
Q_b=f_bP_b+g_bP_q, ~~~ Q_q=f_qP_b+g_qP_q;
\end{array}
\label{PQ}
\end{equation}
 here  $\lambda=S_x$. We say that system (\ref{main}) possesses a  `generic' Lax pair if relations (\ref{PQ}) are in involution, and $P, Q$ do not satisfy any extra first-order relations other than (\ref{PQ}). Generic Lax pair is a generic solution to (\ref{PQ}). Differentiating each of these relations by $a, b, p, q$ and $\lambda$, we obtain 30 relations which, in the non-degenerate case, can be uniquely resolved for all second-order partial derivatives of $P$ and $Q$, thus leading to a closed system. It can be verified directly that the resulting system  is involutive if and only if the functions $f$ and $g$ satisfy  integrability conditions (\ref{*}).
This finishes the proof of Theorem \ref{thm2}.

\medskip

\noindent {\bf Remark 1}. The general solution of the involutive system for $P$ and $Q$ depends on $6$ arbitrary constants. This comes from the count  $2+10-6=6$, where $2$ corresponds to the values of  $P$ and $Q$, $10$ comes from the values of their first-order derivatives, and $-6$ is due to the $6$ first-order relations (\ref{PQ}).  On the other hand, the  Lax pair is form-invariant under transformations
$$
S\to \alpha S+\beta u+\gamma v+\mu x+\nu y+\eta t,
$$
which also involve $6$ arbitrary constants. Modulo these transformations,   `generic'  Lax pair is essentially unique.

\medskip

\noindent {\bf Remark 2}. One can show that, as a consequence of relations (\ref{PQ}), the quantities $P_{\lambda}$ and $Q_{\lambda}$ satisfy the dispersion relation,
\begin{equation}
\det \left[ \left(\begin{array}{cc}
f_a & f_p\\
g_a & g_p
\end{array}\right)
+
\left(\begin{array}{cc}
f_b & f_q\\
g_b & g_q
\end{array}\right)P_{\lambda}
-
\left(\begin{array}{cc}
1 & 0\\
0 & 1
\end{array}\right)Q_{\lambda}
\right]
=0.
\label{det}
\end{equation}
This relation plays an important role in the geometric interpretation of Lax pairs discussed below.

\subsection{Lax pairs and null totally geodesic surfaces}
\label{sec:null}

Here we outline a general construction that leads from  Lax pair (\ref{Lax1}) to null totally geodesic  surfaces of the corresponding Einstein-Weyl structure $(\mathbb{D}, \ g, \ \omega)$. Differentiating (\ref{Lax1}) by $x$ and setting $S_x=\lambda$ we obtain
\begin{equation}
\lambda_y=P_{\lambda} \lambda_x+P_{a}a_x+P_{b}b_x+P_pp_x+P_qq_x, ~~~
\lambda_t=Q_{\lambda} \lambda_x+Q_{a}a_x+Q_{b}b_x+Q_pp_x+Q_qq_x.
\label{lambda}
\end{equation}
With this system we associate two  vector fields,
$$
\begin{array}{c}
X=\frac{\partial}{\partial y}-P_{\lambda}\frac{\partial}{\partial x}+(P_{a}a_x+P_{b}b_x+P_pp_x+P_qq_x)\frac{\partial}{\partial \lambda},
\  \\
Y=\frac{\partial}{\partial t}-Q_{\lambda}\frac{\partial}{\partial x}+(Q_{a}a_x+Q_{b}b_x+Q_pp_x+Q_qq_x)\frac{\partial}{\partial \lambda},
\end{array}
$$
which live in the extended four-dimensional space with coordinates $x, y, t, \lambda$. Note that the compatibility condition, $\lambda_{yt}=\lambda_{ty}$, is equivalent to the commutativity of these vector fields: $[X, Y]=0$.
The geometry behind this construction is as follows. Let us consider the cotangent bundle $Z$ of a solution $u(x, y, t), \ v(x, y, t)$, with local coordinates
$(x, y, t, S_x, S_y, S_t)$. Equations (\ref{Lax1}) specify a four-dimensional submanifold $M^4\subset Z$ parametrised by $x, y, t$ and $ \lambda$. The compatibility of  equations (\ref{Lax1}) indicates that this submanifold is coisotropic.
Vector fields $X, Y$  generate  the kernel of the restriction to $M^4$ of the symplectic form $dS_x\wedge dx+dS_y\wedge dy+dS_t\wedge dt$. Equations (\ref{lambda}) mean that $X, Y$ are tangential to the Lagrangian submanifold of $M^4$ defined by the equation $\lambda=\lambda(x, y, t)$.  Solutions $S(x, y, t)$ of  equations (\ref{Lax1}) can be interpreted as generating functions of Lagrangian submanifolds contained in $M^4$. Such submanifolds depend on one arbitrary function of a single variable; they can be obtained by taking any one-parameter subfamily of the two-parameter family of integral surfaces of the distribution $\langle X, Y\rangle$ foliating $M^4$.

Projecting the two-parameter family of integral surfaces of the distribution $\langle X, Y\rangle$ from $M^4$  to the space of independent variables $x, y, t$ we obtain a two-parameter family of  null totally geodesic  surfaces of the Weyl connection $\mathbb{D}$. Indeed,  projecting $X$ and $Y$ gives vector fields
$$
\hat X=\frac{\partial}{\partial y}-P_{\lambda}\frac{\partial}{\partial x}, ~~~
\hat Y=\frac{\partial}{\partial t}-Q_{\lambda}\frac{\partial}{\partial x},
$$
which commute if and only if $\lambda$ satisfies equations (\ref{lambda}). It remains to show  that $\langle\hat X, \hat Y\rangle $ is a null distribution (that is,  tangential to the null cones of $g$), and  that the covariant derivatives  $\mathbb{D}_{\hat X}\hat X, \ \mathbb{D}_{\hat X}\hat Y,\ \mathbb{D}_{\hat Y} \hat X, \ \mathbb{D}_{\hat Y}\hat Y $ belong to $\langle \hat X, \hat Y\rangle$.
Equivalently, one can introduce the covector $\theta=dx+P_{\lambda}dy+ Q_{\lambda}dt$ that annihilates $\hat X, \hat Y$, and verify that  $\theta$ is null, and that $\mathbb{D}_{\hat X}\theta \wedge \theta=\mathbb{D}_{\hat Y}\theta \wedge \theta=0$.
This follows from  equations (\ref{PQ}) satisfied by the functions $P(\lambda, a, b, c)$ and $Q(\lambda, a, b, c)$. In particular, the fact that $\theta$ is null follows from  identity (\ref{det}).



\section{${\bf SL}(5)$ as the equivalence group:  geometry of  ${\bf Gr}(3, 5)$}
\label{sec:equiv}

Let  $V^5$ be 5-dimensional vector space with  coordinates $x^1, x^2, x^3, u, v$,  set ${\bf p}=(u, v)^t$ and  ${\bf x}=(x^1, x^2, x^3)^t$.  Three-dimensional linear subspaces of $V$  are defined by the equation $d{\bf p}=U d{\bf x}$ where
$$
U=\left( \begin{array}{ccc}
 u_1 & u_2 &u_3
 \\
v_1 & v_2 & v_3
 \end{array}\right).
$$
 Thus, the  Grassmannian ${\bf Gr}(3, 5)$ is six-dimensional, with an affine chart  identified with the  space of $2\times 3$ matrices $U$. 
 Solutions to system (\ref{main}),  $u(x^1, x^2, x^3)$ and $v(x^1, x^2, x^3)$, can be interpreted as three-dimensional  submanifolds of $V^5$.  Their tangent spaces  are specified by matrices $U$ whose entries coincide with first-order derivatives of $u$ and $v$: $u_i=\partial u /\partial x^i, \ v_i=\partial v/\partial x^i$.
Thus, system (\ref{main}) defines a fourfold $X\subset  {\bf Gr}(3, 5)$, while solutions to (\ref{main}) correspond to three-dimensional submanifolds of $V^5$ whose Gaussian images are contained in $X$.
The  action of  ${\bf SL}(5)$,
$$
\left(\begin{array}{c}
d{\bf \tilde p}
 \\
 d{\bf \tilde x}
 \end{array}\right)=
\left( \begin{array}{cc}
 A & B
 \\
C & D
 \end{array}\right)
\left( \begin{array}{c}
d{\bf  p}
 \\
 d{\bf  x}
 \end{array}\right),
 $$
 naturally extends to ${\bf Gr}(3, 5)$:
 \begin{equation}
 \tilde U=(AU+B)(CU+D)^{-1};
 \label{Sp}
 \end{equation}
here $A, B, C, D$ are $2\times 2, \ 2\times 3, \ 3\times 2$ and $3\times 3$ matrices, respectively; notice that  the extended action is no longer linear. Transformation law (\ref{Sp}) suggests that the action of ${\bf SL}(5)$ preserves the class of equations (\ref{main}), indeed, first-order derivatives transform through  first-order derivatives only. Transformations of this form preserve  integrability, so that the group ${\bf SL}(5)$ can be viewed as the natural {\it equivalence group} of the problem: it maps integrable systems to integrable systems. Thus,  ${\bf SL}(5)$ is the point symmetry group of  integrability conditions (\ref{*}) derived in Sect. \ref{sec:hydro}. The classification of integrable systems  (\ref{main}) will be performed modulo this equivalence: two ${\bf SL}(5)$-related equations should be regarded as `the same'.

Geometrically, our problem is reduced to the classification of fourfolds $X$ of the  Grassmannian ${\bf Gr}(3, 5)$ that  satisfy certain `integrability' conditions  (to be specified later, see Sect. \ref{sec:secant}), modulo the action of ${\bf SL}(5)$.


Given two infinitesimally close three-dimensional linear subspaces  of $V^5$ defined by $2\times 3$ matrices $U$ and $U+dU$, the condition that their   intersection is  two-dimensional is given by ${\it rank}~  dU=1$. This condition is  invariant under the action of ${\bf SL}(5)$ as specified by (\ref{Sp}). Thus, each tangent space of ${\bf Gr}(3, 5)$ is equipped with the  cone $C$ defined by the equation ${\it rank}~ dU=1$, that is, by the vanishing of all $2\times2$ minors of $dU$. Projectivisation of this  cone is known as the Segre variety; it  is a non-singular algebraic threefold of degree three. The field of Segre cones supplies ${\bf Gr}(3, 5)$ with the generalised flat conformal structure; it  is manifestly invariant under the action of ${\bf SL}(5)$. The converse is also true: transformations from ${\bf SL}(5)$ are the {\it only}  diffeomorphisms  that preserve the  field of Segre cones.

\medskip

\noindent {\bf Proposition 1.} {\it The group of conformal automorphisms of the field of  Segre cones is isomorphic to ${\bf SL}(5)$.}

\medskip

\noindent  This is a well-known fact:  direct proof would consist of the calculation of conformal automorphisms of the  family of Segre cones.  Let us  point out that, in coordinates $u_i, v_i$, the infinitesimal generators corresponding to equivalence transformations (\ref{Sp}) are as follows:

\noindent {\it 6 translations:}
\begin{eqnarray*}
&&{\bf U}_{i}=\frac{\partial}{\partial{u_{i}}},  ~~~ {\bf V}_{i}=\frac{\partial}{\partial{v_{i}}} ,
\end{eqnarray*}

\noindent {\it 12 linear transformations (note the relation $\sum {\bf X}_{ii}={\bf L}_{11}+{\bf L}_{22}$):}
\begin{eqnarray*}
&&{\bf X}_{ij}= u_{i}\frac{\partial}{\partial{u_{j}}}+v_{i}\frac{\partial}{\partial{v_{j}}},  ~~
{\bf L}_{11}= u_{k}\frac{\partial}{\partial{u_{k}}}, ~~ {\bf L}_{12}=u_{k}\frac{\partial}{\partial{v_{k}}}, ~~ {\bf L}_{21}= v_{k}\frac{\partial}{\partial{u_{k}}}, ~~ {\bf L}_{22}= v_{k}\frac{\partial}{\partial{v_{k}}}.
\end{eqnarray*}

\noindent {\it 6 non-linear (projective) transformations:}
\begin{eqnarray*}
&&{\bf P}_i=u_iu_k\frac{\partial}{\partial{u_{k}}}+v_iu_k\frac{\partial}{\partial{v_{k}}}, ~~
{\bf Q}_i=u_iv_k\frac{\partial}{\partial{u_{k}}}+v_iv_k\frac{\partial}{\partial{v_{k}}},
\end{eqnarray*}

\medskip

\noindent The ideal $\Omega$ defining the field of Segre cones is generated by  quadratic forms,
$$
\Omega={\rm span} ~ \{ du_idv_j-dv_idu_j\},
$$
which are nothing but second fundamental forms of the Pl\"ucker embedding of the Grassmannian, see e.g. \cite{Griffiths}. It remains to point out that any vector field $X$ satisfying $L_X\Omega=0 \ ({\rm mod}\  \Omega$) is spanned by the above infinitesimal generators.
We  refer to \cite{Goncharov, Bertram, Gindikin, Gindikin1} for  generalisations of this Liouville-type result.

Consider the action of the equivalence group ${\bf SL}(5)$ on the space $J^1(\mathbb{R}^4,\mathbb{R}^2)$
of 1-jets of functions $f,g$ of  variables $a,b,p,q$.
This is a 14-dimensional space with coordinates $u_i, v_i, f, g, f_{u_i}, f_{v_i}$, $g_{u_i}, g_{v_i}$, $i=1,2$, which
can be viewed as  an affine chart in the bundle  of 4-dimensional
tangent subspaces of  ${\bf Gr}(3, 5)$. The action of ${\bf SL}(5)$ is  canonically defined
in the latter space, but by abuse of notation we will be working with $J^1$.

\medskip

\noindent {\bf Lemma.} {\it The group ${\bf SL}(5)$ has a unique Zariski open orbit in 
the space $J^1(\mathbb{R}^4,\mathbb{R}^2)$
(its complement consists of 1-jets of degenerate systems).}

\medskip

The proof of this  statement is as follows. The group ${\bf SL}(5)$ acts transitively on 
the Grassmannian ${\bf Gr}(3, 5)$
with the stabilizer of a point $o$ being the parabolic subgroup 
$P_o=S({\bf GL}(2)\times {\bf GL}(3))\ltimes(\mathbb{R}^2\otimes\mathbb{R}^3)$ of upper-triangular block matrices 
of the size $2+ 3$. As only $S({\bf GL}(2)\times {\bf GL}(3))$ acts on $T_o{\bf Gr}(3, 5)$, 
(i.e. $\mathbb{R}^2\otimes\mathbb{R}^3$ acts trivially),  the action is transitive on
 4-planes corresponding to non-degenerate 1-jets of $f, g$.
At the level of Lie algebra $\mathfrak{sl}(5)$,  prolongation of the above infinitesimal generators
 to $J^1(\mathbb{R}^4,\mathbb{R}^2)$ has full rank in the Zariski open set of non-degenerate 
1-jets. Indeed, the $24\times 14$  matrix of coefficients of these vector fields drops  rank precisely 
on the submanifold $\det g^{ij}=0$ where $g^{ij}$ is the conformal structure defined by the characteristic variety. 

\medskip

\noindent  {\bf Remark 1.} 
The above lemma allows one to assume that all sporadic factors depending on first-order derivatives of $f$ and $g$ that arise in the process of Gaussian elimination in the proofs of our main results, are nonzero. This considerably simplifies the arguments by eliminating unessential branching.
Furthermore, in the verification of polynomial identities involving first and second-order partial derivatives of $f$ and $g$ one can, without any loss of generality, give the first-order derivatives any `generic' numerical values: this often renders otherwise impossible computations manageable. 

\medskip

In fact, even the action of ${\bf SL}(3)\subset P_o$ on $T_o{\bf Gr}(3, 5)$ is already transitive, 
and  the stabilizer of a non-degenerate four-plane in $T_o{\bf Gr}(3, 5)$ is the subgroup 
$H_o={\bf GL}(2)\ltimes(\mathbb{R}^2\otimes\mathbb{R}^3)$ of $P_o$,
where we identify $\mathbb{R}^3\simeq S^2(\mathbb{R}^2)$ as an ${\bf SL}(2)$-representation. Prolongation of this 
action to $J^2(\mathbb{R}^4,\mathbb{R}^2)$ is locally free near a generic point (straightforward calculation shows that the prolongation of infinitesimal generators of $\mathfrak{sl}(5)$ to the space of 2-jets has full rank 24 at generic points).

\medskip

\noindent {\bf Proposition 2.} {\it The action of the equivalence group ${\bf SL}(5)$ on the $30$-dimensional parameter space of integrable systems (\ref{main}) is locally free (in the Zariski open set of generic points). }

\medskip

That is, `generic' orbits of this action are $24$-dimensional, so that  `generic' integrable systems of type (\ref{main}) do not possess any continuous point symmetries coming from the equivalence group. Thus, the moduli space of integrable systems depends on $30-24=6$ essential parameters. To prove the above statement let us note that the submanifold in the space of 3-jets of $f, g$ given by integrability conditions (\ref{*}) projects diffeomorphically onto the space of 2-jets $J^2(\mathbb{R}^4,\mathbb{R}^2)$ (which has dimension 30). 


\medskip

\com{Given system (\ref{evol}),
$$
u_3=f(u_1, u_2, v_1, v_2), ~~~ v_3=g(u_1, u_2, v_1, v_2),
$$
the Lie algebra of its equivalence group results from the above infinitesimal generators upon substitution $u_3=f, \ v_3=g$. Prolonging these generators to the first-order jet space $J^1$, one obtains 24 vector fields on the 14-dimensional space with coordinates $u_i, v_i, f, g, f_{u_i}, f_{v_i}, g_{u_i}, g_{v_i}$, here $i=1, 2$. One can verify that, for non-degenerate systems, these vector fields span $J^1$. To be precise, the $24\times 14$ matrix formed by coefficients of these vector fields has rank 14 almost everywhere; the rank drops  on the submanifold $\det g^{ij}=0$ where $g^{ij}$ denotes  principal symbol. The complement to this submanifold in $J^1$ is the open orbit of ${\bf SL}(5)$-action. This allows one to assume that, modulo  equivalence transformations, all sporadic factors depending on first-order derivatives of $f$ and $g$ that arise in the process of Gaussian elimination in the proofs of our main results, are nonzero.
This considerably simplifies the proofs by eliminating unessential branching. Similarly, in the verification of polynomial identities involving first- and second-order partial derivatives of $f$ and $g$, without any loss of generality one can give the first-order derivatives any `generic' numerical values: this often makes otherwise impossible computations manageable.
}

\noindent {\bf Remark 2.} The action of ${\bf SL}(5)$ on the 30-dimensional parameter space $M^{30}$  is algebraic,
so that by the Rosenlicht theorem \cite{Rosen}  there exists a quotient  $M^{30}_{\op{reg}}/{\bf SL}(5)$, which is a rational
algebraic variety of dimension 6 (the geometric quotient $M^{30}/{\bf SL}(5)$ is singular). 
Indeed, the equivalence group transforms the initial conditions of system (19) rationally.

\com{By the arguments
before Proposition 2, it coincides with the quotient of the fiber of the projection 
$J^2(\mathbb{R}^4,\mathbb{R}^2)\to J^1(\mathbb{R}^4,\mathbb{R}^2)$ over the point $o$ by the affine  
action of  $H_o$.}

\section{Examples and classification results}
\label{sec:ex}

This section contains various classification results based on  integrability conditions (\ref{*}).


\noindent
We produce an  abundance of non-trivial examples of integrable systems of type (\ref{main}), both known and new, expressible in elementary functions, theta functions,  modular forms and generalised hypergeometric functions.

\subsection{Linearisable systems}
\label{sec:lin}

In this section we  characterise  systems (\ref{main}) that can be linearised by a transformation from the equivalence group ${\bf SL}(5)$. Taking a linear system, say  $u_3=v_1, \ v_3=u_2$, and applying transformations from the equivalence group, one obtains  systems of  Monge-Amp\`ere type,
\begin{equation}
\begin{array}{c}
a^{ij}(u_iv_j-u_jv_i)+b^iu_i+c^iv_i+m=0, \\
\alpha^{ij}(u_iv_j-u_jv_i)+\beta^iu_i+\gamma^iv_i+\mu=0,
\end{array}
\label{Monge}
\end{equation}
where each equation represents a linear combination of  minors of the $2\times 3$ matrix
$$
U=\left( \begin{array}{ccc}
 u_1 & u_2 &u_3
 \\
v_1 & v_2 & v_3
 \end{array}\right).
$$
Conversely, in three dimensions every system of this form  is linearisable.

\medskip

 \noindent {\bf Proposition 3.} {\it For a non-degenerate system of type (\ref{main}),  the
 following conditions are equivalent:

 \noindent (a) System is linearisable by a transformation from the equivalence group ${\bf SL}(5)$.

  \noindent (b) System belongs to  Monge-Amp\`ere class  (\ref{Monge}).

 \noindent (c) System is invariant under  an $8$-dimensional subgroup of ${\bf SL}(5)$.

 \noindent (d) The characteristic variety  defines conformal structure that is  flat on every solution.
}

\medskip

 \centerline{\bf Proof:}

 \medskip

 \noindent {\bf equivalence $(a)\Longleftrightarrow (b)$}: One needs to show that taking a non-degenerate  linear system, say  $u_3=v_1, \ v_3=u_2$ (all non-degenerate linear systems of type (\ref{main}) are ${\bf SL}(5)$-equivalent), and applying transformations from the equivalence group, one obtains {\it all} systems of Monge-Amp\`ere type. The easiest way to see this is the following. First of all, any non-degenerate linear system is invariant under an $8$-dimensional subgroup of ${\bf SL}(5)$; note that 8 is the maximal possible value for the dimension of the stabiliser of the action of the equivalence group on the space of fourfolds $X\subset\mathbf{Gr}(3,5)$. The corresponding infinitesimal generators can be calculated using the standard machinery of group analysis. Thus, the linear system $u_3=v_1, \ v_3=u_2$ is invariant under the subgroup with 8 infinitesimal generators,
$$
\begin{array}{c}
{\bf U}_{1}, ~~ {\bf V}_{2}, ~~ {\bf U}_{2}+{\bf V}_{3}, ~~ {\bf U}_{3}+{\bf V}_{1}, \\
{\bf X}_{13}+2{\bf X}_{32}+{\bf L}_{12}, ~~ {\bf X}_{23}+2{\bf X}_{31}+{\bf L}_{21}, ~~ {\bf X}_{11}+{\bf X}_{22}+{\bf X}_{33}, ~~
{\bf X}_{11}-{\bf X}_{22}+{\bf L}_{11},
\end{array}
$$
 (we use the notation of Sect. \ref{sec:equiv}). The Lie algebra is isomorphic  to a semi-direct product ${V}_3\rtimes {\mathfrak{gl}}(2)$ where ${V}_3\simeq S^3(\mathbb{R}^2)$ is the irreducible 4-dimensional
representation of $\mathfrak{sl}(2)$ (and hence of $\mathfrak{gl}(2)=\mathfrak{sl}(2)\oplus\mathbb{R}$, with non-trivial action of $\mathbb{R}$).
 Thus, applying to a linear system transformations from the equivalence group, one obtains a variety of systems depending on $24-8=16$ essential parameters. It remains to point out that the class of Monge-Amp\`ere systems also depends on $16$ essential parameters.

\medskip

\noindent {\bf equivalence $(a)\Longleftrightarrow (c)$}: The first implication follows from the fact that any non-degenerate  linear system is invariant under an $8$-dimensional subalgebra of  ${\bf SL}(5)$. To establish the converse, let $G$
be the symmetry group of system (\ref{main}). We can always assume that
the point $o$, specified by $u_{i}=v_i=0$, belongs to the fourfold $X\subset {\bf Gr}(3, 5)$
corresponding to our system. Let $G_o$ be the stabiliser of this
point in $G$. Note that $\dim G - \dim G_o \le4$, as $G$ takes $X$ to itself. The
stabiliser $P$ of the point $o$ is spanned by  infinitesimal generators
${\bf X}_{ij}, \ {\bf L}_{ij}, \ {\bf P}_i, \ {\bf Q}_i$.
Since the system is non-degenerate, we can always bring it to a
canonical form:
\begin{equation}\label{cf}
u_3=v_1+o(u_i, v_i),  ~~~  v_3=u_2+o(u_i, v_i).
\end{equation}
This form (together with the point $o$) is stabilised by the following four elements of $P$:
$$
{\bf X}_{13}+2{\bf X}_{32}+{\bf L}_{12}, ~~ {\bf X}_{23}+2{\bf X}_{31}+{\bf L}_{21}, ~~ {\bf X}_{11}+{\bf X}_{22}+{\bf X}_{33}, ~~
{\bf X}_{11}-{\bf X}_{22}+{\bf L}_{11}.
$$
Thus, $\dim G_o \le4$ so that $\dim G  \le8$. The equality holds only if
$\dim G_o=4$. However, the generator  ${\bf X}_{11}+{\bf X}_{22}+{\bf X}_{33}$ acts by non-trivial
rescalings on  terms of  order 2 and higher in \eqref{cf}. Hence, for $\dim G_o=4$, all higher-order terms must vanish identically, leading to a linear system.

\medskip

\noindent {\bf equivalence $(a)\Longleftrightarrow (d)$}: Representing system (\ref{main}) in evolutionary form (\ref{evol})
and introducing the corresponding conformal structure $g$ (see the proof of Theorem \ref{thm1} for explicit formulae),
the condition responsible for conformal flatness in three dimensions is the vanishing of the Cotton tensor,
\begin{equation}
\nabla_r(R_{pq}-\tfrac{1}{4}Rg_{pq})=\nabla_q(R_{pr}-\tfrac{1}{4}Rg_{pr}),
\label{cot}
\end{equation}
where $R_{pq}$ is the Ricci tensor, $R$ is the scalar curvature, and $\nabla$ denotes covariant differentiation in the Levi-Civita connection of   $g$. Calculating (\ref{cot}) and using  (\ref{evol}) and its differential consequences to eliminate all higher-order partial derivatives of $u$ and $v$ containing differentiation by $x^3$, we obtain expressions that have to vanish identically
in the remaining higher-order derivatives  (no more than fourth-order derivatives  will occur in this calculation).  In particular, equating to zero coefficients at fourth-order derivatives of $u$ and $v$, we obtain the system  of second-order PDEs for $f$ and $g$:
\begin{equation}
\begin{array}{c}
f_{u_iu_i}=\frac{2g_{u_i}}{g_{v_i}-f_{u_i}}f_{u_iv_i}, ~~~ f_{v_iv_i}=\frac{2f_{v_i}}{f_{u_i}-g_{v_i}}f_{u_iv_i}, \\
\ \\
f_{u_iu_j}=\frac{g_{u_j}}{g_{v_i}-f_{u_i}}f_{u_iv_i}+\frac{g_{u_i}}{g_{v_j}-f_{u_j}}f_{u_jv_j}, ~~~
f_{v_iv_j}=\frac{f_{v_j}}{f_{u_i}-g_{v_i}}f_{u_iv_i}+\frac{f_{v_i}}{f_{u_j}-g_{v_j}}f_{u_jv_j}, \\
\ \\
f_{u_iv_j}+f_{u_jv_i}=\frac{f_{u_j}-g_{v_j}}{f_{u_i}-g_{v_i}}f_{u_iv_i}+
\frac{f_{u_i}-g_{v_i}}{f_{u_j}-g_{v_j}}f_{u_jv_j};
\end{array}
\label{fij}
\end{equation}
here $i, j$ take any values from 1 to 2; the equations for $g$ can be obtained by the simultaneous substitution $f\leftrightarrow g$ and $u\leftrightarrow v$. The general solution of  system (\ref{fij}) leads to Monge-Amp\`ere systems. Furthermore, modulo (\ref{fij}) all remaining components of the Cotton tensor  vanish identically (in fact, this follows from the linearisability of Monge-Amp\`ere systems in 3D).
This finishes the proof of Proposition 3.

\medskip

\noindent {\bf Remark.} Equations of Monge-Amp\`ere type  have a clear geometric interpretation.
Recall that the  Grassmannian ${\bf Gr}(3, 5)$ is  (locally) identified  with the space of $2\times 3$  matrices $U$. Minors of $U$ define the Pl\"ucker embedding of ${\bf Gr}(3, 5)$ into projective space $\mathbb{P}^{9}$. We identify ${\bf Gr}(3, 5)$ with the image of this  embedding, which  is a non-singular  algebraic variety of degree five. Thus,   Monge-Amp\`ere systems correspond  to  sections of ${\bf Gr}(3, 5)$ by subspaces $\mathbb{P}^7\subset\mathbb{P}^9$.

\medskip

We  emphasize that the linearisability of Monge-Amp\`ere systems  is an essentially three-dimensional phenomenon: in higher dimensions there exist integrable non-linearisable  examples of Monge-Amp\`ere type, e.g.
$$
u_2-v_1=0, ~~~ u_3v_4-u_4v_3-1=0,
$$
which is equivalent to the Pleba\'nski first heavenly equation \cite{Plebanski}.

\subsection{Linearly degenerate systems}
\label{sec:lindeg}

From the point of view  of their `complexity', linearly degenerate systems come next after linearisable systems. 
In this section we obtain  a complete  list of  normal forms of linearly degenerate integrable systems: all of them come from quadratic maps $\mathbb{P}^4 \dashrightarrow {\bf Gr}(3, 5)$ given by the classical construction of Chasles.

The definition of linear degeneracy is inductive: a multi-dimensional system is said to be linearly degenerate (completely exceptional  \cite{Boillat1}) if  all its traveling wave reductions to two dimensions are linearly degenerate.  Thus, it is sufficient to define this concept in 2D case,
\begin{equation}
u_2=f(u_1, v_1), ~~~ v_2=g(u_1, v_1).
\label{1+1}
\end{equation}
Setting $u_{1}=a, \ v_{1}=p$  and differentiating by $x^1$ one can rewrite this system   in  two-component quasilinear form,
$$
a_2=f(a, p)_1, ~~~ p_2=g(a, p)_1,
$$
or, in matrix notation,
$$
\left(\begin{array}{c}
a\\
p
\end{array}\right)_2=A\left(\begin{array}{c}
a\\
p
\end{array}\right)_1,  ~~~ A=\left(\begin{array}{cc}
f_a & f_p\\
g_a & g_p
\end{array}\right).
$$
Recall that the matrix $A$ is said to be linearly degenerate if its eigenvalues (assumed real and distinct) are constant in the direction of the corresponding  eigenvectors. Explicitly, $L_{r^i}\lambda^i=0$, no summation, where $L_{r^i}$ denotes Lie derivative in the direction of the  eigenvector $r^i$, and $A r^i=\lambda^i r^i$. For quasilinear systems, the property of linear degeneracy is known to be related, under appropriate  `small norm' assumptions, to the impossibility of breakdown of smooth initial data, leading to  global solvability of a generic Cauchy problem \cite{Roz}.
In terms of the original functions $f(u_1, v_1)$ and $g(u_1, v_1)$, the conditions of linear degeneracy reduce to a pair of second-order differential constraints,
\begin{equation}
\begin{array}{c}
(f_{u_1}-g_{v_1})f_{u_1u_1}+2g_{u_1}f_{u_1v_1}+g_{u_1}g_{v_1v_1}+f_{v_1}g_{u_1u_1}=0,\\
\ \\
(g_{v_1}-f_{u_1})g_{v_1v_1}+2f_{v_1}g_{u_1v_1}+f_{v_1}f_{u_1u_1}+g_{u_1}f_{v_1v_1}=0.
\end{array}
\label{lind}
\end{equation}
Equations (\ref{1+1})  define a 2-parameter family of 2-dimensional linear subspaces of a vector space $V^4$ (on projectivisation, a congruence of lines in $\mathbb{P}^3$). Geometrically, equations (\ref{lind}) mean that focal surfaces of this congruence degenerate into curves \cite{Agafonov}. We point out that  system (\ref{lind}) is invariant under the  natural action of the 2D equivalence group ${\bf SL}(4)$. Requiring that all traveling wave reductions of a multi-dimensional system to 2D are linearly degenerate in the above sense, we obtain differential characterisation of linear degeneracy:

\medskip

 \noindent {\bf Proposition 4.} {\it A $d$-dimensional  system,
 $
 u_d=f(u_i, v_i),\ v_d=g(u_i, v_i),
 $
 $i=1, \dots, d-1$, is linearly degenerate if and only if the functions $f$ and $g$ satisfy the relations
 \begin{equation}
\begin{array}{c}
{Sym}_{\{i, j, k\}}\left((f_{u_k}-g_{v_k})f_{u_iu_j}+g_{u_k}(f_{u_iv_j}+f_{u_jv_i})+f_{v_k}g_{u_iu_j}+g_{u_k}g_{v_iv_j}\right)=0,\\
\ \\
{Sym}_{\{i, j, k\}}\left((g_{v_k}-f_{u_k})g_{v_iv_j}+f_{v_k}(g_{u_iv_j}+g_{u_jv_i})+g_{u_k}f_{v_iv_j}+f_{v_k}f_{u_iu_j}\right)=0,
\end{array}
\label{sym}
\end{equation}
where Sym denotes complete symmetrisation over $i, j, k\in\{1, \dots, d-1\}$.

}

\medskip

 \centerline{\bf Proof:}

 \medskip

 \noindent
Let us begin with the $3$-dimensional case,
$$
u_3=f(u_1,  u_2, v_1, v_2), ~~~ v_3=g(u_1, u_2, v_1, v_2).
$$
Looking for a traveling wave reduction in the form
$$
u(x^1, x^2, x^3)=u(y^1, y^2)+y^3, ~~~ v(x^1, x^2, x^3)=v(y^1, y^2)+y^4,
$$
where $y^i$ are arbitrary constant-coefficient linear forms in the independent variables $x^k$, one obtains a system of  type (\ref{1+1}) for
$u(y^1, y^2)$ and $v(y^1, y^2)$. It is required to be linearly degenerate for {\it every} choice of linear forms $y^i$. This condition imposes a set of  second-order differential constraints for  $f$ and $g$,
\begin{equation}
\begin{array}{c}
(f_{u_j}-g_{v_j})f_{u_iu_i}+2(f_{u_i}-g_{v_i})f_{u_iu_j}+2g_{u_j}f_{u_iv_i}+2g_{u_i}(f_{u_iv_j}+f_{u_jv_i})+\\
\ \\
f_{v_j}g_{u_iu_i}+2f_{v_i}g_{u_iu_j}+g_{u_j}g_{v_iv_i}+2g_{u_i}g_{v_iv_j}=0,\\
\ \\
(g_{v_j}-f_{u_j})g_{v_iv_i}+2(g_{v_i}-f_{u_i})g_{v_iv_j}+2f_{v_j}g_{u_iv_i}+2f_{v_i}(g_{u_iv_j}+g_{u_jv_i})+\\
\ \\
g_{u_j}f_{v_iv_i}+2g_{u_i}f_{v_iv_j}+f_{v_j}f_{u_iu_i}+2f_{v_i}f_{u_iu_j}=0,
\end{array}
\label{lind2}
\end{equation}
(no summation), where  $i, j\in \{1, 2\}$ (8 relations altogether). Notice that the second half of these relations follows from the first one under the simultaneous substitution $f\leftrightarrow g, \ u\leftrightarrow v$.
Relations (\ref{lind2}) can be obtained from (\ref{sym}) by setting   $k=j$. For $i=j=1$ relations (\ref{lind2}) simplify to (\ref{lind}).
Relations (\ref{lind2}) can be derived  in either of the two  ways:

\medskip

\noindent {\bf Method 1.}  Applying  the full 3D equivalence group ${\bf SL}(5)$ to  constraints (\ref{lind}) one gets 8 linearly independent relations that are equivalent to (\ref{lind2}). One can show that the set of relations (\ref{lind2}) is  ${\bf SL}(5)$-invariant. This procedure generalises to higher dimensions in an obvious way. For instance, in 4D,  applying transformations from ${\bf SL}(6)$ to  constraints (\ref{lind}) one gets 20 linearly independent relations  (\ref{sym}), $d=4$.

\medskip

\noindent {\bf Method 2.}  Looking for particular traveling wave reductions of the form
$u(x^1, x^2, x^3)=u(y^1, y^2)$,  $v(x^1, x^2, x^3)=v(y^1, y^2)$ where $y^1=x^1+\lambda x^2, \ y^2=x^3$, one obtains the reduced system  in the form
$$
u_{y^2}=F(u_{y^1}, v_{y^1})=f(u_{y^1},  \lambda u_{y^1}, v_{y^1}, \lambda v_{y^1}), ~~~ v_{y^2}=G(u_{y^1}, v_{y^1})=g(u_{y^1},  \lambda u_{y^1}, v_{y^1}, \lambda v_{y^1}).
$$
Here $F$ and $G$ have to satisfy  relations (\ref{lind}). Since $F_{u_{y^1}}=f_{u_1}+\lambda f_{u_2}, $ $F_{u_{y^1}u_{y^1}}=f_{u_1 u_1}+2\lambda f_{u_1 u_2}+\lambda^2 f_{u_2 u_2}, $ etc,
both relations (\ref{lind}) will be  polynomial of degree 3 in $\lambda$. Equating to zero coefficients of these polynomials we get all  relations (\ref{lind2}). This  generalises to higher dimensions: thus,
in 4D, given a system
$$
u_4=f(u_1,  u_2, u_3, v_1, v_2, v_3), ~~~ v_4=g(u_1, u_2, u_3, v_1, v_2, v_3),
$$
one looks for traveling wave reductions
$u(x^1, x^2, x^3, x^4)=u(y^1, y^2),$  $v(x^1, x^2, x^3, x^4)=v(y^1, y^2)$ where $y^1=x^1+\lambda x^2+\mu x^3, \ y^2=x^4$. For the reduced system, both relations (\ref{lind}) become polynomials of degree 3 in $\lambda$ and $\mu$, each having $10$ coefficients. Equating them to zero we obtain $20$ differential constraints (\ref{sym}) constituting  conditions of linear degeneracy in 4D. This finishes the proof of Proposition 4.

\medskip
\noindent {\bf Remark 1.} Second-order relations (\ref{lind2}) governing linearly degenerate systems in 3D are not in involution, and their prolongation implies all third-order integrability conditions (\ref{*}). This requires differentiating equations (\ref{lind2}) three times and solving for the fifth, fourth and third-order partial derivatives of $f$ and $g$ in terms of the first and second-order derivatives. 

 \medskip

 \noindent {\bf Proposition 5.} {\it The parameter space of linearly degenerate integrable systems in 3D is $22$-dimensional. Furthermore, the following conditions are equivalent:

\noindent {(a)} System (\ref{main}) is linearly degenerate and integrable.

\noindent  {(b)} There exists a unique flat symmetric connection on $X$ in which the  associated ${\bf GL}(2, \bbbr)$ structure is parallel.

}

\medskip

 \centerline{\bf Proof:}

 \medskip

 \noindent One can verify that 8 constraints (\ref{lind2}) are compatible with  integrability conditions (\ref{*}). Remembering that the parameter space of integrable systems (\ref{main})  is $30$-dimensional, we thus obtain a $30-8=22$-dimensional parameter space of 3D linearly degenerate integrable systems.

The equivalence of (a) and (b) can be demonstrated as follows. Using parametric equations of the fourfold $X$ as in Theorem \ref{thm1},
$$
u_{x}=a,\ u_{y}=b, \ v_{x}=p,\ v_{y}=q, \ u_{t}= f(a,b,  p, q), \ v_t=g(a, b, p, q),
$$
we define cubic cones of the induced ${\bf GL}(2, \bbbr)$ structure as the intersection of three quadratic forms, $\omega^{\alpha}=0$, where
$$
\omega^1=dadq-dbdp, ~~~ \omega^2=dadg-dpdf, ~~~ \omega^3=dbdg-dqdf.
$$
Note that these quadratic forms are nothing but restrictions to $X$ of the second fundamental forms of the Grassmannian
${\bf Gr}(3, 5)\subset \mathbb{P}^9$. A connection $\nabla$ preserving this ${\bf GL}(2, \bbbr)$ structure is defined by the relations $\nabla \omega^{\alpha}=0 ~ {\rm mod}  \langle \omega^{\beta}\rangle$. Under the assumption that $\nabla$ is symmetric, these relations lead to a linear inhomogeneous system for the 40 Christoffel's symbols of $\nabla$. This system is uniquely solvable if and only if the functions $f$ and $g$ satisfy  conditions of linear degeneracy (\ref{lind2}). One can verify by direct calculation that the vanishing of the curvature of this symmetric connection is equivalent to  integrability conditions (\ref{*}).
 This finishes the proof of Proposition 5.

\medskip

Our next goal is to provide a complete list of normal forms of linearly degenerate integrable systems. The key example generating an open part of the  22-dimensional parameter space is as follows:

\medskip

\noindent {\bf Example.} Let $a_1, a_2, a_3$ and  $\tilde a_1, \tilde a_2, \tilde a_3$ be  two triplets of numbers  such that $a_1+a_2+a_3=0$ and $\tilde a_1+\tilde a_2+\tilde a_3$=0. One can show that the following system,
\begin{equation}
a_1
\tilde a_2u_xv_y- a_2\tilde a_1u_yv_x=0,
~~~ a_1\tilde a_3 u_xv_t-a_3\tilde a_1u_t v_x=0, \label{aa}
\end{equation}
 is linearly degenerate and integrable. Eliminating $v$ we obtain a second-order equation for $u$,    $a_1u_xu_{yt}+a_2u_yu_{xt}+a_3u_tu_{xy}=0$.
Similarly, eliminating $u$ we obtain the analogous equation for $v$,  $\tilde a_1v_xv_{yt}+\tilde a_2v_yv_{xt}+\tilde a_3v_tv_{xy}=0$. This construction first appeared in \cite{Zakharevich} in the context of bi-Hamiltonian systems and Veronese webs in 3D, see also \cite{DK}.  Rewriting  system (\ref{aa}) in  simplified form,
$$
{u_x}v_y=\alpha {u_y}{v_x}, ~~~  {u_xv_t}=\beta {u_tv_x},
$$
$\alpha, \beta = const$, one can show that:

\noindent (i) Systems with generic values of $\alpha$ and $\beta$ are not  ${\bf SL}(5)$-equivalent.

\noindent (ii) For  generic choice of $\alpha$ and $\beta$ the corresponding system is invariant under a $4$-dimensional subgroup of the equivalence group (with infinitesimal generators ${\bf X}_{11}$,  ${\bf X}_{22}$,  ${\bf X}_{33}$,  ${\bf L}_{11}$ in the notation of Sect. \ref{sec:equiv}, that form a Cartan subalgebra of ${\bf SL}(5)$), and hence generates a $24-4=20$-dimensional orbit.

\noindent Acting on $(\alpha, \beta)$-systems  by transformations from the equivalence group one thus obtains a $22$-parameter family generating an open part of the parameter space of linearly degenerate integrable systems.

\medskip

All  linearly degenerate integrable systems in 3D, including the above example, can be obtained from the following geometric construction.
Consider  projective space $\mathbb{P}^4$ with homogeneous coordinates $\xi=(\xi^1:\xi^2:\xi^3:\xi^4:\xi^5)$. Let $A$ be a projective automorphism of $\mathbb {P}^4$ defined by a $5\times 5$ matrix from ${\bf SL}(5)$. Consider the family of lines
$\overline{(\xi,\eta)}$
through $\xi$ and $\eta=\xi A$ (the locus of lines spanned by an argument and the value of a projective transformation was apparently first discussed by Chasles \cite{Chasles}; see also \cite{Dolgachev},  p. 556).  The  Pl\"ucker coordinates $p^{ij}=\xi^i\eta^j-\xi^j\eta^i$  are the $2\times 2$ minors of the matrix
$$
\left(\begin{array}{ccccc}
\xi^1 & \xi^2 & \xi^3 & \xi^4 & \xi^5\\
\eta^1 & \eta^2 & \eta^3 & \eta^4 & \eta^5
\end{array}\right);
$$
they define quadratic map from the family of lines
$\overline{(\xi, \eta)}$
(which is itself isomorphic to $\mathbb {P}^4$) into the Grassmannian of lines in $\mathbb {P}^4$, that is, into ${\bf Gr}(2, 5)$.
 By duality, this gives a fourfold $X$ in ${\bf Gr}(3, 5)$, and the corresponding  system (\ref{main}). 
 Fourfolds $X$ arising from this construction are images of quadratic maps  $\mathbb{P}^4 \dashrightarrow {\bf Gr}(3, 5)\subset \mathbb{P}^9$. Explicit parametric equations of  $X$, as well as of the corresponding system (\ref{main}), can be obtained from the factorised representation
$$
\left(\begin{array}{ccccc}
\xi^1 & \xi^2 & \xi^3 & \xi^4 & \xi^5\\
\eta^1 & \eta^2 & \eta^3 & \eta^4 & \eta^5
\end{array}\right)=
\left(\begin{array}{cc}
\xi^1 & \xi^2 \\
\eta^1 & \eta^2
\end{array}\right)
\left(\begin{array}{ccccc}
1 & 0 & u_1 & u_2 & u_3\\
0 & 1 & v_1& v_2 & v_3
\end{array}\right),
$$
explicitly,
 \begin{equation}
u_1=\frac{p^{32}}{p^{12}},\ u_2=\frac{p^{42}}{p^{12}},\ u_3=\frac{p^{52}}{p^{12}}, \
v_1=\frac{p^{13}}{p^{12}},\ v_2=\frac{p^{14}}{p^{12}},\ v_3=\frac{p^{15}}{p^{12}}.
 \label{uvxi}
 \end{equation}
In the (generic) diagonal case, $A=diag(\lambda^i)$, this gives
$$
u_1=\frac{\lambda^2-\lambda^3}{\lambda^2-\lambda^1}\frac{\xi^2\xi^3}{\xi^1\xi^2}, ~~~
u_2=\frac{\lambda^2-\lambda^4}{\lambda^2-\lambda^1}\frac{\xi^2\xi^4}{\xi^1\xi^2}, ~~~
u_3=\frac{\lambda^2-\lambda^5}{\lambda^2-\lambda^1}\frac{\xi^2\xi^5}{\xi^1\xi^2},
 $$
$$
v_1=\frac{\lambda^3-\lambda^1}{\lambda^2-\lambda^1}\frac{\xi^1\xi^3}{\xi^1\xi^2}, ~~~
v_2=\frac{\lambda^4-\lambda^1}{\lambda^2-\lambda^1}\frac{\xi^1\xi^4}{\xi^1\xi^2}, ~~~
v_3=\frac{\lambda^5-\lambda^1}{\lambda^2-\lambda^1}\frac{\xi^1\xi^5}{\xi^1\xi^2},
 $$
leading to the  $(\alpha, \beta)$-system, $u_1v_2=\alpha u_2v_1, ~ u_1v_3=\beta u_3v_1$, with
$$
\alpha =\frac{(\lambda^2-\lambda^4)(\lambda^3-\lambda^1)}{(\lambda^2-\lambda^3)(\lambda^4-\lambda^1)}, ~~~
\beta =\frac{(\lambda^2-\lambda^5)(\lambda^3-\lambda^1)}{(\lambda^2-\lambda^3)(\lambda^5-\lambda^1)}.
$$
The corresponding fourfold $X$ coincides with the quadratic  image  of $\mathbb{P}^4$ defined by a linear system of quadrics that pass through 5 points in general position (in our parametrisation, quadrics $\xi^i\xi^j, \ i\ne j$, that pass through 5 base points each having only one nonzero  coordinate  $\xi$). We point out that smooth quadratic embeddings of projective spaces into Grassmannians were discussed in \cite{ugaglia}. It is important to emphasize that our $X$ is not smooth:  lines through any pair of base points in $\mathbb{P}^4$ correspond to  singular points of  $X$. In algebro-geometric language,
quadratic maps from $\mathbb {P}^n$ to the Grassmannian of lines in $\mathbb{P}^n$ are given by  rank 2 vector bundles E with the first Chern class equal to 2, together with n+1 sections that generate E. In the Chasles construction,
$E = {\cal O}(1)+{\cal O}(1)$, and the space of sections is the graph of the map $A: H^0({\cal O}(1))\to H^0({\cal O}(1))$.

\medskip

Table 1 below comprises a complete list of 7 canonical forms of linearly degenerate integrable systems (\ref{main}) labelled by Segre types (essentially, Jordan normal forms) of $5\times 5$ matrices $A$ from the Chasles construction. Thus, type $[3, 2]$ denotes  two Jordan blocks of  size $3\times 3$ and $2\times 2$, type $[5]$ - one $5\times 5$ Jordan block, etc (to get canonical forms as presented in Table 1 we take $A$ in  Jordan normal form, eliminate $\xi$'s from (\ref{uvxi}), and use equivalence transformations to remove unessential parameters). The last column provides dimensions of  stabilisers of the resulting  systems under the action of the equivalence group
${\bf SL}(5)$.

\medskip

\begin{center}
 \centerline{\footnotesize{Table 1: Canonical forms of linearly degenerate integrable systems  in 3D}}
 \medskip
 \begin{tabular}{ | l | l | l | p{1.1cm} |} \hline
 {\footnotesize  Segre type}  & {\footnotesize Canonical form}  & {\footnotesize Dim(stab)} \\
 \hline

{\footnotesize  [11111]}&\footnotesize{$u_1v_2=\alpha u_2v_1$} & ~~~~ 4  \\
&{\footnotesize $u_1v_3=\beta u_3v_1$} & \\
\hline

{\footnotesize  [2111]}&\footnotesize{$u_1v_2-u_2v_1=v_1v_2$} & ~~~~ 4  \\
&{\footnotesize $u_1v_3-u_3v_1=\alpha v_1v_3$} & \\
\hline

{\footnotesize  [221]}&\footnotesize{$u_1v_2-u_2v_1=v_1^2$} & ~~~~ 4  \\
&{\footnotesize $u_1v_3-u_3v_1=v_1v_3$} & \\
\hline

{\footnotesize  [311]}&\footnotesize{$u_2=v_1v_2$} & ~~~~ 4  \\
&{\footnotesize $u_3=(1-v_1)v_3$} & \\
\hline

{\footnotesize  [32]}&\footnotesize{$u_2=v_1v_2$} & ~~~~ 5  \\
&{\footnotesize $u_3=v_2+v_1v_3$} & \\
\hline

{\footnotesize  [41]}&\footnotesize{$u_1=v_2-v_1^2$} & ~~~~ 5  \\
&{\footnotesize $u_3=(1-v_1)v_3$} & \\
\hline

{\footnotesize  [5]}&\footnotesize{$u_1=v_2-v_1^2$} & ~~~~ 6  \\
&{\footnotesize $u_2=v_3-v_1v_2$} & \\
\hline
\end{tabular}

\end{center}
Systems presented in Table 1 are not ${\bf SL}(5)$-equivalent. Even though all 4-dimensional stabilisers are commutative subalgebras of ${\bf SL}(5)$, they are not conjugate (in particular, only the first of them is a Cartan subalgebra).

\medskip

\noindent{\bf Remark 2.} Note that Table 1 does not contain a  linear system. It can be recovered from the following  generalisation of the Chasles construction.
Take  projective space $\mathbb{P}^4$ with homogeneous coordinates $\xi=(\xi^1:\xi^2:\xi^3:\xi^4:\xi^5)$,  let $A$ and $B$ be two   $5\times 5$ matrices from ${\bf SL}(5)$. Consider the family of lines
$\overline{(\eta,\varphi)}$
through $\eta=\xi A$ and $\varphi=\xi B$. A linear system corresponds to the case  when the pencil $A+\lambda B$ consist of two Kronecker blocks of the size $1\times 2$ and $4\times 3$ (Kronecker blocks of other size lead to degenerate systems).

\medskip

\noindent{\bf Proposition 6.} {\it
The Chasles construction gives all linearly degenerate integrable systems in 3D.
}

 \medskip

\centerline{\bf Proof:}

\medskip

Consider the Grassmann variety of all pencils $A+\lambda B$ in the space of $5\times 5$ matrices. Denote by $\Sigma$ the open subset in ${\bf Gr}(2,25)$ which corresponds to non-degenerate equations via the Chasles construction. Clearly, this is an irreducible algebraic variety.
The Chasles construction defines a rational algebraic map from ${\bf Gr}(2,25)$ to the 22-dimensional variety of all linearly degenerate integrable systems, with $\Sigma$ as the set of regular points.  The fibres of this map are   orbits of the right action of  ${\bf SL}(5)$ on pencils $A+\lambda B$. Moreover, we can easily verify that  this action has a discrete stabiliser in the case when the pencil $A+\lambda B$ is of Segre type $[11111]$.
Therefore,  generic fibres of this map are $24$-dimensional. Thus, its image is an irreducible variety of dimension $2\cdot 23 - 24 = 22$. As ${\bf Gr}(2,25)$ is a projective variety, the image has to coincide with the set of all linearly degenerate integrable systems.

\medskip

Table 2 brings together some further examples of linearly degenerate integrable systems that can be interpreted as B\"acklund transformations: on elimination of $v$, they lead to a second-order PDE for $u$, similarly, on elimination of $u$ they lead to a second-order PDE for $v$ (columns 2 and 3).  Each of these systems is equivalent to one of the canonical forms presented in Table 1 (as indicated by its type).

\medskip

\begin{center}
 \centerline{\footnotesize{Table 2: Examples of B\"acklund transformations in 3D}}
 \medskip
 \begin{tabular}{ | l | l | l | p{1.1cm} |} \hline
 {\footnotesize  System (\ref{main})}  & {\footnotesize Equation for $u$}  & {\footnotesize Equation for $v$} \\ 
 $$ &  &  \\
  \hline

{\footnotesize Type\ [11111]}& &   \\
{\footnotesize $ a_1\tilde a_3v_tu_x-a_3\tilde a_1 v_xu_t=0$}& {\footnotesize $a_1u_xu_{yt}+a_2u_yu_{xt}+a_3u_tu_{xy}=0$ }&{\footnotesize $\tilde a_1v_xv_{yt}+\tilde a_2v_yv_{xt}+\tilde a_3v_tv_{xy}=0$ }    \\
{\footnotesize $a_1\tilde a_2v_yu_x- a_2\tilde a_1v_xu_y=0$ }&  &   \\ \hline

{\footnotesize Type\ [2111]}& &   \\
 {\footnotesize $ (\lambda-1)v_y-u_yv_x=0$}&{\footnotesize $u_{yt}+u_yu_{xt}-u_tu_{xy}=0$} &{\footnotesize $v_xv_{yt}+(\lambda-1)v_yv_{xt}-\lambda v_t v_{xy}=0$} \\ 
{\footnotesize $\lambda v_t-u_tv_x=0$} &  &   \\ \hline

{\footnotesize Type\ [221]}& &   \\
{\footnotesize $v_y-\lambda u_yv_x=0$}& { \footnotesize $u_{yy}+u_yu_{xt}-u_tu_{xy}=0$} & {\footnotesize $v_tv_{xy}-v_yv_{xt}+\lambda(v_yv_{xy}-v_xv_{yy})=0$ }   \\
{\footnotesize $v_t+(\lambda^2 u_y-\lambda u_t)v_x=0$ }&  &   \\ \hline

{\footnotesize Type\ [221]}& &   \\
 {\footnotesize $u_yv_t-1=0$}& {\footnotesize $u_{xy}+u_yu_{xt}-u_xu_{yt}=0$} & {\footnotesize $v_{xt}+v_tv_{xy}-v_x v_{yt}=0$}     \\
 {\footnotesize $u_x-u_yv_x=0$} &  &   \\ \hline

{\footnotesize Type\ [311]}& &   \\
 {\footnotesize $v_t+(\lambda -u_y)v_y=0$}& {\footnotesize $u_{xt}+u_xu_{yy}-u_yu_{xy}=0$ }& {\footnotesize $v_tv_{xy}-v_yv_{xt}+\lambda(v_yv_{xy}-v_xv_{yy})=0$ }    \\
 {\footnotesize $\lambda v_x-u_xv_y=0$ }&  &   \\ \hline

{\footnotesize Type\ [32]}& &   \\
 {\footnotesize $u_yv_x-1=0$}& {\footnotesize $u_{yy}+u_yu_{xt}-u_tu_{xy}=0$} & {\footnotesize $v_{xt}+v_xv_{yy}-v_y v_{xy}=0$}   \\
 {\footnotesize $u_t-u_yv_y=0$} &  &   \\ \hline

{\footnotesize Type\ [41]}& &   \\
 {\footnotesize $v_y+(\lambda -u_x)v_x=0$}& {\footnotesize $u_{xt}+u_{yy}+u_yu_{xx}-u_xu_{xy}=0$ }& {\footnotesize $\left(\frac{v_y}{v_x}\right)_y=\left(\frac{\lambda v_y-v_t}{v_x}\right)_x$}     \\
 {\footnotesize $v_t+(\lambda^2-\lambda u_x+u_y)v_x=0$}&  &   \\ \hline

\end{tabular}

\end{center}
Some of the second-order PDEs from Table 2 have appeared in different contexts in \cite{Pavlov, Shabat1, Zakharevich, Bur, Odesskii3, Kras}.

\medskip
\noindent {\bf Remark 3.}  The equivalence group acts algebraically on the parameter space of linearly generate integrable systems, so we can consider the rational quotient $M^{22}_\text{reg}/\mathbf{SL}(5)$. Since generic orbits of this action have dimension 20, the moduli space of linearly degenerate integrable systems is two-dimensional.

\subsection{Further examples  in terms of  modular forms and theta functions}
\label{sec:mod}

This section contains a list of more exotic integrable systems that are not expressible in elementary functions. These examples demonstrate that one should not expect any `simple  parametrisation' of integrable systems of type ({\ref{main}).

\medskip

\noindent{\bf Example 1.} Let us begin with the system
$$
u_t=v_x, ~~~ v_{t}=\frac{u_{y}}{v_{x}}+\frac{1}{6}\eta (u_{x})v_{x}^{2},
$$
which appeared in  \cite{MaksEgor} in the classification of integrable hydrodynamic chains. In this case the integrability conditions result in the Chazy equation for $\eta$,  
$$
\eta^{\prime \prime \prime }+2\eta \eta ^{\prime \prime }-3{\eta^{\prime} }^2=0,
$$
whose generic solution is known to be  the weight two Eisenstein series $E_2$ associated with the full modular group.

\medskip

\noindent {\bf Example 2.} Consider  the system
$$
v_{x}+u_{x}u_{y}r(u_{t})=0, ~~~ v_t=u_y.
$$
In this case the integrability conditions result in the  third-order ODE for $r$,
\begin{equation*}
r'''(r'-r^2) -r''^2+ 4{r}^3{r''}+2r'^3 - 6{r}^2{r'}^2= 0,
\end{equation*}
which  appeared recently in  different context  in the theory of  modular forms of level two: compare with  equation (4.7) from \cite{A2}. Its  generic solution is given by the  Eisenstein series
$$
r(u_t)=1-8\sum_{n=1}^{\infty}\frac{(-1)^nne^{4nu_t}}{1-e^{4nu_t}}, 
$$
 associated with the congruence subgroup $\Gamma_0(2)$ of the modular group.

\medskip

\noindent {\bf Example 3.} As a generalisation of Example 2, let us consider the system
$$
v_{x} + u_{x} r(u_{y}, u_{t})= 0, ~~~ v_t=u_y.
$$
In a somewhat different representation, it first appeared in \cite{Bur}.  It was demonstrated that the integrability conditions imply  $r=2s'/s$ where prime denotes differentiation by $u_{y}$, and
$ s(u_{y},u_{t}) = \theta\left(\frac{u_{y}}{2\pi},-\frac{u_{t}}{\pi i} \right)$
 is the Jacobi theta function:
\begin{equation*} \theta(z,\tau) = 1 + 2\sum_{n=1}^\infty e^{\pi i n^2 \tau} \cos(2\pi n z).
\end{equation*}

\medskip

\noindent {\bf Example 4.} Further generalisation,
$$
v_{x}+f(u_{x}, u_{y}, u_{t})=0, ~~~ v_t=u_y,
$$
 was discussed in \cite{Fer4},  where it was shown that the requirement of integrability implies that generic $f$ is given by the ratio of two Jacobi theta functions:
$$
f(u_{x}, u_{y}, u_{t})=-\frac{1}{4}\ln \frac{\theta_1(u_{t}, u_{y}-u_{x})}{\theta_1(u_{t}, u_{y}+u_{x})}.
$$

\subsection{Integrable systems in terms of generalised hypergeometric functions}
\label{sec:hyper}

In this section we  utilise the  construction of \cite{Odesskii2, Odesskii3} which parametrises   integrable quasilinear equations in $2+1$ dimensions in terms of solutions of the generalised hypergeometric system,
\begin{equation} \label{hyper}
\begin{array}{c}
\displaystyle \frac{\partial^2 h}{\partial z_i \partial
z_j}=\frac{s_i}{z_i-z_j} \frac{\partial h}{\partial
z_j}+\frac{s_j}{z_j-z_i} \frac{\partial h}{\partial
z_i},  ~~~  i\ne j, \\
\ \\
\displaystyle \frac{\partial^2 h}{\partial z_i^2}=-\Big(1+\sum_{j=1}^{n+2} s_j\Big) \frac{s_i}{z_i (z_i-1) }\, h+
\frac{s_i}{z_i (z_i-1)} \sum_{j\ne i}^n \frac{z_j
(z_j-1)}{z_j-z_i}\cdot
\frac{\partial h}{\partial z_j}+\\[7mm]
\displaystyle \Big(\sum_{j\ne i}^n \frac{s_j}{z_i-z_j}+
\frac{s_i+s_{n+1}}{z_i}+ \frac{s_i+s_{n+2}}{z_i-1}\Big)
\frac{\partial h}{\partial z_i},
\end{array}
\end{equation}
here  $h$ is a function of $n$ variables $z_1, \dots, z_n$, and $s_1,...,s_{n+2}$ are arbitrary constants. This system is involutive and possesses $n+1$ linearly independent solutions  known as  generalised hypergeometric functions \cite{gel, Odesskii2}.
Here we consider the case  $n=4$. Let  $\{h^1, h^2, g^1, g^2, g^3\}$ be a basis of solutions of  system (\ref{hyper}). Let us introduce the  parametric  formulae
\begin{equation}
\begin{array}{c}
\displaystyle u_x=\frac{W(g^1,  h^1)}{W(h^1, h^2)}, ~~~ u_y=\frac{W(g^2, h^1)}{W(h^1, h^2)}, ~~~ u_t=\frac{W(g^3, h^1)}{W(h^1, h^2)},\\
\ \\
\displaystyle v_x=\frac{W(g^1,  h^2)}{W(h^1,  h^2)}, ~~~ v_y=\frac{W(g^2, h^2)}{W(h^1, h^2)}, ~~~ v_t=\frac{W(g^3, h^2)}{W(h^1, h^2)},
\end{array}
\label{par}
\end{equation}
where $W$ is the Wronskian,  $W(f, g)=fg'-gf'$, and prime denotes differentiation with respect to one of the variables $z_i$ (for definiteness, we  assume $'=\partial/\partial {z_4}$). Formulae (\ref{par}) can be viewed as parametric equations of a fourfold  $X\subset {\bf Gr}(3, 5)$. Eliminating the parameters $z_1, \dots, z_4$ we obtain  system of type (\ref{main}).

\medskip

\noindent {\bf Proposition 7.} Parametric equations (\ref{par}) define  {\it generic integrable} system of type (\ref{main}).

\medskip

\centerline{\bf Proof:}

\medskip

Given  hypergeometric system ({\ref{hyper}), the  construction of  \cite{Odesskii2} requires an additional ingredient, namely, the choice of a $k$-dimensional subspace of solutions of  system (\ref{hyper}).
We will only need a particular case of the general scheme  that corresponds to  $n=4, \ k=2$. Thus, let $h^1, h^2$ be a pair of linearly independent solutions of  system (\ref{hyper}). Let us complete this pair to a basis  $\{h^1, h^2, g^1, g^2, g^3\}$, and introduce the following system
of 6 first-order quasilinear PDEs for  $z_1, \dots, z_4$,  considered as functions of the auxiliary  variables $x, y, t$:
$$
\begin{array}{c}
\left(\frac{W(g^1,  h^1)}{W(h^1,  h^2)}\right)_y=\left(\frac{W(g^2, h^1)}{W(h^1, h^2)}\right)_x, ~~~
\left(\frac{W(g^1,  h^1)}{W(h^1,  h^2)}\right)_t=\left(\frac{W(g^3, h^1)}{W(h^1, h^2)}\right)_x, ~~~
\left(\frac{W(g^2, h^1)}{W(h^1, h^2)}\right)_t=\left(\frac{W(g^3, h^1)}{W(h^1, h^2)}\right)_y, \\
\ \\
\left(\frac{W(g^1,  h^2)}{W(h^1,  h^2)}\right)_y=\left(\frac{W(g^2, h^2)}{W(h^1, h^2)}\right)_x, ~~~
\left(\frac{W(g^1,  h^2)}{W(h^1,  h^2)}\right)_t=\left(\frac{W(g^3, h^2)}{W(h^1, h^2)}\right)_x, ~~~
\left(\frac{W(g^2,  h^2)}{W(h^1,  h^2)}\right)_t=\left(\frac{W(g^3, h^2)}{W(h^1, h^2)}\right)_y.
\end{array}
$$
It was demonstrated in \cite{Odesskii2} that this system is integrable by the method of hydrodynamic reductions, and possesses a dispersionless Lax pair. The conservative structure of this system implies the existence of potentials $u$ and $v$ specified by (\ref{par}). Thus, parametric equations (\ref{par}) indeed  give rise to  {\it integrable} system (\ref{main}).  Since 
the moduli space of integrable systems (\ref{main}) is 6-dimensional and, for $n=4$,  hypergeometric system (\ref{hyper}) depends on 6 essential parameters $s_1, \dots, s_6$, the claim follows.
One can show that transformations from the equivalence group ${\bf SL}(5)$ are in one-to-one correspondence with  linear transformations of  the  basis $\{h^1, h^2, g^1, g^2, g^3\}$.

We  emphasise that, although formulae (\ref{par}) parametrise   {\it generic} integrable systems (\ref{main}), the degeneration procedure leading to examples from Section \ref{sec:mod} is far from trivial.


\section{Geometric aspects of integrability}
\label{sec:geom}

In this section we adopt a geometric point of view and consider 3D system (\ref{main}) as equations defining a fourfold $X$ in the   Grassmannian ${\bf Gr}(3, 5)$. Our aim  is to reformulate the integrability conditions in intrinsic geometric terms.

\medskip

\noindent Recall that  in Sect. \ref{sec:equiv} we introduced the main object needed to develop  differential geometry  of submanifolds of ${\bf Gr}(3, 5)$, namely,   the field of  Segre cones in the tangent bundle   T${\bf  Gr}(3, 5)$. The field of Serge cones  is invariant under the natural action of the  equivalence group ${\bf SL}(5)$, and plays the role of an ambient  flat conformal structure.

\medskip

\noindent In Sect. \ref{sec:gl2} we develop  geometry of fourfolds  $X\subset {\bf Gr}(3, 5)$. The induced field of Segre cones gives rise to  a field of twisted cubics  in the projectivised  tangent bundle  $\mathbb{P}{\rm T}X$. Thus, a fourfold in  ${\bf Gr}(3, 5)$ carries an intrinsic  ${\bf GL}(2, \bbbr)$ structure.

\medskip

\noindent In Sect. \ref{sec:secant} we show that two-component and three-component hydrodynamic reductions of system (\ref{main}) correspond to bisecant  and holonomic trisecant submanifolds of the fourfold $X$, respectively. This allows us to reformulate  the integrability geometrically as the existence of an infinity of holonomic trisecant
submanifolds parametrised by three arbitrary functions of one variable.

\medskip

\noindent In Sect. \ref{sec:confsymp} we demonstrate that every fourfold $X$ corresponding to  integrable system  (\ref{main}) carries a canonical conformally symplectic structure.

\medskip

\noindent In Sect. \ref{sec:curvtor} we reformulate the integrability conditions in terms of the torsion/curvature of the Bryant connection of the associated ${\bf GL}(2, \bbbr)$ structure.

\subsection{${\bf GL}(2, \bbbr)$ structure on a fourfold $X\subset {\bf Gr}(3, 5)$}
\label{sec:gl2}

Let $X$ be a non-degenerate fourfold in the  Grassmannian ${\bf Gr}(3, 5)$. Taking a point $o\in X$ and projectivising the intersection of the tangent space ${\rm T}_oX$ with the Serge cone $C$ in ${\rm T}_o{\bf Gr}(3, 5)$  one obtains a twisted cubic, that is, a rational normal curve of degree three. Indeed, the projectivised Serge cone is a non-singular rational cubic threefold in $\mathbb{P}^5$, so that its section by the projectivised tangent space $\mathbb{P}{\rm T}_oX=\mathbb{P}^3$ is a non-singular curve of degree three, that is, a twisted cubic. This twisted cubic can also be interpreted as the set of matrices of rank one in the tangent space ${\rm T}_oX$, when we identify ${\rm T}_o{\bf Gr}(3, 5)$ with the space of $2\times 3$  matrices.
Thus, the projectivised tangent bundle  of the fourfold $X$ is equipped with a field of  twisted cubics.  This supplies $X$ with a ${\bf GL}(2, \bbbr)$ structure.

Alternatively, one can say that each tangent space to $X$ is  identified with a four-dimensional space of binary cubics: in this picture the rational normal curve corresponds to cubics with a triple root.
In representation-theoretic language, the tangent space $T_oX$ is an irreducible $\mathfrak{sl}(2)$
module $V_3$; here and below $V_l\simeq S^l(\mathbb{R}^2)$ is the irreducible $\slg(2)$-representation
of dimension $l+1$.

\medskip

\noindent {\bf Remark 1.} Another natural source of  ${\bf GL}(2, \bbbr)$ structures in  four dimensions  is provided by  4-th order ODEs  with vanishing W\"unschmann invariants \cite{Bryant, Doubrov, Dun, Nur}. In this context, ${\bf GL}(2, \bbbr)$ structures are induced on the  parameter spaces of solutions to the corresponding ODEs. We emphasise that the structures coming from ODEs always have zero torsion (so that all geometry is contained in  the curvature), while in our case  the torsion is generally nonzero, and defines the curvature.

\noindent {\bf Remark 2.} An interesting  class of  ${\bf GL}(2, \bbbr)$ structures in  five  dimensions appeared recently in the context of  integrable equations of the dispersionless Hirota type. These structures are  induced on hypersurfaces of the Lagrangian Grassmannian that correspond to such equations, see \cite{Hirota, Smith} for further details.

\medskip

The most important property of  ${\bf GL}(2, \bbbr)$ structures in  four dimensions is the existence of  a unique intrinsic  connection (the Bryant connection) whose torsion takes values in the 8-dimensional irreducible representation of ${\bf GL}(2, \bbbr)$ \cite{Bryant}.
Let us recall the construction.
The cubic cones of the  ${\bf GL}(2, \bbbr)$ structure induced on $X$ are defined as the intersection of 3 quadrics, $\omega^{\alpha}=0$, where
 $$
\omega^1=dadq-dbdp, ~~~ \omega^2=dadg-dpdf, ~~~ \omega^3=dbdg-dqdf;
 $$
here $f(a, b, p, q)$ and $g(a, b, p, q)$ correspond to the right-hand sides of  system (\ref{evol}).
The tangent space $T_oX$ is an irreducible $\mathfrak{sl}(2)$ module $V_3$, and since
$\Lambda^2V_3=V_0+V_4$, it has a canonical conformal almost symplectic structure $\Omega$
(see Sect. \ref{sec:confsymp} for explicit formulae).
Raising one index of  $\omega^{\alpha}$ by means of $\Omega$ we obtain 3 operators
$A^{\alpha}\in \gl(4),\ \alpha=1, 2, 3,$ that form a basis  of the irreducible embedding of $\slg(2)$
into $\gl(4)$. Define the bilinear form
 \[
B^{\alpha \beta} = \tr A^{\alpha} A^{\beta},
 \]
which, up to a constant factor, is the Killing form of $\slg(2)$. Let $B^\#$ be the inverse of $B$, and let
 \[
C = 20 B^\#_{\alpha \beta} A^{\alpha} A^{\beta}
 \]
be the (normalized) Casimir operator. Then $C$ is a well-defined element of the universal enveloping algebra of $\slg(2)$ which is independent of the choice of the basis $\omega^{\alpha}$ (the coefficient $20$ appears here because we use the embedding of $\slg(2)$ into $\gl(4)$, rather than  the adjoint representation, to define the Casimir operator).
The Casimir $C$ acts naturally on any irreducible $\slg(2)$ module $V_l$ as a scalar operator (due to the Schur lemma), 
\[
C_{V_l} = \lambda_l\operatorname{Id},\quad \lambda_l=l(l+2).
\]
The Bryant connection $\nabla$ is uniquely defined by the following two properties \cite{Bryant}:
\begin{itemize}
\item $\nabla \omega^{\alpha}=0 ~ {\rm mod}\ \langle \omega^{\beta}\rangle $. This means that $\nabla$ preserves the  ${\bf GL}(2, \bbbr)$ structure.

\item The torsion $T$ of the Bryant connection $\nabla$  lies in the unique 8-dimensional submodule $V_7$ of the $\slg(2)$-representation in the space of all algebraic torsions. This can be expressed via the following linear condition:
\[
C\cdot T = \lambda_7 T=63 T.
\]
\end{itemize}

In more detail, the first condition gives
$56$ equations for the $64$ Christoffel's symbols $\Gamma^i_{jk}$ of the connection $\nabla$. Of these $48$ are linearly independent, so that $\nabla$ exists with the freedom of $64-48=16$ arbitrary
functions.
To fix these functions we note that since the tangent space $\tau=T_oX$ is an irreducible $\mathfrak{sl}(2)$ module $V_3$,  the torsion $T$ of $\nabla$ decomposes as follows:
 $$
\Lambda^2\tau^*\otimes\tau=(V_0\oplus V_4)\otimes V_3=V_1\oplus 2V_3\oplus V_5\oplus V_7.
 $$
We can change the connection by a $\mathfrak{gl}(2)$-equivariant gauge from ${\rm Hom}(\tau, \mathfrak{gl}(2))$,
which decomposes as $V_3\otimes(V_0+V_2)=V_1\oplus 2V_3\oplus V_5$. This removes the  freedom of 16 arbitrary functions and yields a unique connection with torsion $T\in V_7$.  Note that the
action of the operators $A^{\alpha}$ on tensors is standard, for instance, for the torsion tensor one has
$(A^\alpha\cdot T)^k_{ij}=(A^\alpha)^k_aT^a_{ij}-(A^\alpha)^a_iT^k_{aj}-(A^\alpha)^a_jT^k_{ia}.$ This defines the action of the universal enveloping algebra, in particular, of the Casimir operator $C$.

\subsection{Geometric interpretation of the integrability conditions}
\label{sec:secant}

As explained in Sect. \ref{sec:gl2}, each projectivised tangent space ${\rm T}_oX$ carries a rational normal curve $\gamma$ of degree three.

\medskip

\noindent {\bf Definition.}

\noindent (a)  {\it Bisecant surface} is a two-dimensional submanifold $\Sigma^2\subset X$ whose projectivised tangent planes are bisecant lines of $\gamma$.

\noindent (b) Given a three-dimensional submanifold $\Sigma^3\subset X$, each tangent space ${\rm T}_o\Sigma^3$ carries three distinguished directions, namely those corresponding to the three points of intersection of $\mathbb{P}{\rm T}_o\Sigma^3$ with $\gamma$.
These directions define a net on $\Sigma^3$, we will call it the {\it characteristic} net.   {\it Holonomic} submanifolds $\Sigma^3$  are defined by the requirement that   the characteristic net is  holonomic (that is, locally a coordinate net).


\medskip

\noindent {\bf Proposition 8.} {\it Bisecant surfaces and holonomic submanifolds  $\Sigma^3$   correspond to two- and three-component hydrodynamic reductions of  system (\ref{main}). Furthermore,

\noindent {\bf (i)} every fourfold $X$  possesses infinitely many bisecant surfaces parametrised  by two arbitrary functions of one variable;

\noindent{\bf (ii)} a fourfold $X$ corresponds to an integrable system  if and only if it possesses infinitely many holonomic submanifolds $\Sigma^3$  parametrised by three arbitrary functions of one variable. Thus, the existence of holonomic submanifolds $\Sigma^3$  is a geometric characterisation of integrability.}

\medskip

\centerline {\bf Proof:}

\medskip

\noindent We follow the notation of Sect. \ref{sec:hydro}. Let us represent our system  in  quasilinear form (\ref{quasi}), and take an $N$-component reduction specified by
$$
a=a(R^{1},..., R^{N}), ~~ b=b(R^{1},...,R^{N}), ~~ p=p(R^{1},...,R^{N}), ~~ q=q(R^{1},...,R^{N}).
$$
 The geometric image  in the  Grassmannian ${\bf Gr}(3, 5)$ is a submanifold $\Sigma^N\subset X$ (it would be sufficient for our purposes to restrict to   $N=2, 3$), represented by the symmetric matrix
\begin{equation}
U=
\left(\begin{array}{ccc}
a & b & f \\
p & q & g
\end{array}
\right),
\label{U}
\end{equation}
parametrised by $R^1, ..., R^N$. Using  equations (\ref{E1}), (\ref{E2})  one obtains the following expression for the derivative of $U$ with respect to $R^i$,
\begin{equation}
\partial_iU=
\left(\begin{array}{ccc}
\partial_ia & \mu^i\partial_ia & \lambda^i \partial_ia\\
\partial_ip & \mu^i\partial_ip & \lambda^i \partial_ip \\
\end{array}
\right).
\label{Ui}
\end{equation}
Thus, $\partial_iU$ is a matrix of rank one, so that the projectivisation of $\partial_iU$ belongs to the rational normal curve $\gamma$, and  coordinates $R^i$ provide the characteristic net on $\Sigma^N$. For $N=2$ we have a two-dimensional surface $\Sigma^2\subset X$ parametrised by $R^1, R^2$. Since both $\partial_1U$ and $\partial_2U$ have rank one,   the surface $\Sigma^2$ is bisecant. As explained in  Sect. \ref{sec:hydro}, every system of type (\ref{main}) (not necessarily integrable) possesses infinitely many two-component reductions parametrised by two arbitrary functions of one variable. Thus, every fourfold $X$ of the  Grassmannian ${\bf Gr}(3, 5)$ possesses infinitely many bisecant surfaces. This establishes  part {(i)} of the Proposition.

Similarly, every three-component reduction corresponds to a holonomic submanifold $\Sigma^3\subset X$ parametrised by  $R^1, R^2, R^3$. Since an integrable system possesses (by definition) infinitely many three-component reductions parametrised by three arbitrary functions of one variable, the corresponding fourfold $X$ possesses infinitely many  holonomic submanifolds $\Sigma^3$. This establishes the first part of (ii).

To finish the proof one needs to show that, conversely,   bisecant surfaces (holonomic submanifolds $\Sigma^3$) of $X$ correspond to two-component (three-component) reductions of the associated system. This can be demonstrated as follows. Let $\Sigma^2$ be a bisecant surface represented in  form (\ref{U}), referred to its characteristic net $R^1, R^2$. Thus,  $a, b,  p, q, f, g$ are functions of $R^1, R^2$ such that the rank of $\partial_iU$ equals one, so that  one can introduce  parametrisation (\ref{Ui}). Compatibility conditions for the equations $\partial_ib=\mu^i\partial_ia$ and $\partial_iq=\mu^i\partial_ip$ imply
\begin{equation*}
\partial_i\partial_ja=\frac{\partial_j\mu^i}{\mu^j-\mu^i}\partial_ia+\frac{\partial_i\mu^j}{\mu^i-\mu^j}\partial_ja ~~~ {\rm and} ~~~
\partial_i\partial_jp=\frac{\partial_j\mu^i}{\mu^j-\mu^i}\partial_ip+\frac{\partial_i\mu^j}{\mu^i-\mu^j}\partial_jp,
\end{equation*}
respectively. Similarly,  compatibility conditions for the equations $\partial_if=\lambda^i\partial_ia$ and $\partial_ig=\lambda^i\partial_ip$ imply
\begin{equation*}
\partial_i\partial_ja=\frac{\partial_j\lambda^i}{\lambda^j-\lambda^i}\partial_ia+\frac{\partial_i\lambda^j}{\lambda^i-\lambda^j}\partial_ja ~~~ {\rm and} ~~~
\partial_i\partial_jp=\frac{\partial_j\lambda^i}{\lambda^j-\lambda^i}\partial_ip+\frac{\partial_i\lambda^j}{\lambda^i-\lambda^j}\partial_jp.
\end{equation*}
Subtracting these equations from each other we get
$$
s_{ij}\partial_ia+s_{ji}\partial_ja=0, ~~~ s_{ij}\partial_ip+s_{ji}\partial_jp=0,
$$
where $s_{ij}=\partial_j\mu^i/(\mu^j-\mu^i)-\partial_j\lambda^i/(\lambda^j-\lambda^i)$. The case $s_{ij}=s_{ji}=0$ gives   commutativity conditions (\ref{comm}), so that we  recover all equations governing two-component  reductions. Thus, bisecant surfaces indeed correspond to two-component reductions.  If one of the coefficients $s_{ij}$ or $s_{ji}$ is nonzero, then the functions $a$ and $p$ must be functionally dependent, so that $dp=c\,da$ for some coefficient $c$. But then we must also have $dq=c\,db$ and $dg=c\,df$, so that the whole tangent space of the surface $\Sigma^2$ consists of rank one matrices. This, however, is not possible for non-degenerate systems (\ref{main}).

For holonomic submanifolds $\Sigma^3$, considerations are essentially the same. Since the existence of three-component reductions is necessary and sufficient for the integrability, this finishes the proof of  Proposition 8.

\subsection{Canonical conformal almost symplectic structure}
\label{sec:confsymp}

 Let $\mathbb{P}^3$ be a projective space with homogeneous coordinates $(x^0 : x^1 : x^2 : x^3)$.  A rational normal curve $\gamma$ of degree three has the standard form $\gamma=(1 : t : t^2 : t^3)$. One can verify that, up to a scalar factor,  there exists a unique $2$-form, $\omega= dx^0\wedge dx^3-3dx^1\wedge dx^2$, for which all tangent planes to the cone over $\gamma$ are Lagrangian, that is,
 $\omega (\gamma, \gamma')=0$. The matrix of $\omega$  is
$$
\omega= \left(
\begin{array}{cccc}
0 & 0 & 0& 1 \\
0 & 0 & -3&0 \\
0 & 3 & 0&0 \\
-1&0&0&0
\end{array}
\right).
$$
Since each projectivised tangent space of $X$ contains a rational normal curve, the `curved' version of the above construction supplies $X$ with a nondegenerate $2$-form $\Omega$  defined up to a conformal factor, that is, a conformal almost symplectic structure.

The differential of $\Omega$ can be represented in the form $d\Omega=\varphi\wedge \Omega$ where $\varphi$ is the so-called Lee form of $\Omega$ (in four dimensions such representation exists, and is unique). Since the rescaling $\Omega \to f\Omega$ gives rise to the transformation  $\varphi \to\varphi+df/f$, the differential $d\varphi$ is independent of the conformal factor $f$. In particular, the condition $d\varphi=0$ is (locally) equivalent to the existence of a function $f$ such that $f\Omega$ is symplectic (such structures $\Omega$ are also known as conformally symplectic).

\medskip

\noindent {\bf Proposition 9.}  {\it The  canonical structure $\Omega$ on the fourfold $X$
corresponding to integrable system (\ref{main}) is locally conformally symplectic. If $H^1(X,\mathbb{R})=0$,
and $X$ is  non-degenerate, then there exists a symplectic structure $\Omega$ on $X$ defined
canonically up to multiplication by a constant. It is parallel with respect to the Bryant connection.}

\medskip

\centerline {\bf Proof:}

\medskip

\noindent Following the notation of Sect. \ref{sec:hydro}  we parametrise   $X\subset {\bf Gr}(3, 5)$ by local coordinates $a, b, p, q$ in the form
$$
X=\left(\begin{array}{ccc}
a&b&f(a, b, p, q)\\
p&q&g(a, b, p, q)\end{array}\right).
$$
The field of  rational normal curves $\gamma$ in the projectivised tangent space of $X$ is specified by the equations ${rk}(dX)=1$, that is
 $$
da\,dq-db\,dp=0, ~~ da\,dg-df\,dp=0, ~~ db\,dg-df\,dq=0,
 $$
or, explicitly,
$
\gamma=\left(m : n:  mt : nt \right)
$
where  $m=g_b+(g_q-f_b)t-f_qt^2, \ n= -g_a+(f_a-g_p)t+f_pt^2$, and $t$ is a parameter.  In matrix notation  we have
$
\gamma=(1:t:t^2:t^3)\cdot A$, where 
$$
A=\left(
\begin{array}{cccc}
g_b & -g_a &0 &0\\
g_q-f_b&f_a-g_p&g_b&-g_a\\
-f_q&f_p&g_q-f_b&f_a-g_p\\
0&0&-f_q&f_p
\end{array}
 \right).
$$
Thus, the matrix of the  $2$-form $\Omega$  is given by $\Omega=A^{-1}\omega (A^{-1})^T$. It remains to calculate the corresponding $1$-form $\varphi$, and verify that the condition $d\varphi=0$  holds identically modulo  integrability conditions (\ref{*}) satisfied by $f$ and $g$.  We point out that $det \Omega=9/ (det A)^2$, furthermore, the condition
$det A=0$ is equivalent to the reducibility of  dispersion relation (\ref{conic}). Thus, our assumption of non-degeneracy of the characteristic variety  is equivalent to non-degeneracy of $\Omega$. Finally, one can verify by direct calculation that the associated symplectic structure, which is unique up to a constant factor, is parallel in the Bryant connection:  this follows from the fact that both the  connection and the symplectic structure are defined in a canonical way.
This finishes the proof of Proposition 9.



\medskip

We emphasize that the condition $d\varphi=0$ is necessary, but not generally sufficient for integrability. For instance, one can show that this condition is satisfied for every system of the form
$u_t=f(u_x)+v_y, ~ u_y=v_x$, where the function $f$ is arbitrary, whereas the integrability implies $f'''=0$. On the other hand, for systems of the form $u_t=v_y, ~ v_t=u_y+g(u_x)$, the condition $d\varphi=0$ implies $g'''g'-2g''^2=0$, which is precisely the integrability condition. We skip details of these calculations.

\subsection{Integrability conditions via differential invariants} 
\label{sec:curvtor}

Let $\nabla$ be the Bryant connection of a ${\bf GL}(2, \bbbr)$-structure on $X$. Its basic invariant is the torsion  $T\in V_7$.
Other tensorial invariants include the curvature $R$, and the covariant derivative
of  torsion, $\nabla T$.
Let us give here explicit coordinate formulae that show our conventions (for non-symmetric connections):
 \begin{gather*}
\nabla_k\partial_i=\Gamma^a_{ik}\partial_a,\quad
\nabla_kdx^i=-\Gamma^i_{ak}dx^a,\
T^k_{ij}=\langle\nabla_i\partial_j-\nabla_j\partial_i,dx^k\rangle=\Gamma^k_{ji}-\Gamma^k_{ij},\\
R^k_{lij}=\langle[\nabla_i,\nabla_j]\partial_l,dx^k\rangle=
\partial_i\Gamma^k_{lj}-\partial_j\Gamma^k_{li}+\Gamma^a_{lj}\Gamma^k_{ai}-\Gamma^a_{li}\Gamma^k_{aj},\\
(\nabla T)^k_{lij}=\partial_lT^k_{ij}+\Gamma^k_{al}T^a_{ij}-\Gamma^a_{il}T^k_{aj}-\Gamma^a_{jl}T^k_{ia}.
 \end{gather*}
Notice that both $R$ and $\nabla T$ (which are third-order differential invariants of the corresponding systems (\ref{evol}), as
 $\Gamma^k_{ij}$ are second-order quantities in partial derivatives of $f$ and $ g$),
belong to the space $\tau\otimes\tau^*\otimes\Lambda^2\tau^*$ which has the following
$\mathfrak{sl}(2)$ decomposition into irreducibles (with multiplicities):
 \begin{equation}\label{dec-tten}
\tau\otimes\tau^*\otimes\Lambda^2\tau^*=2V_0\oplus 4V_2\oplus 5V_4\oplus 4V_6\oplus 2V_8\oplus V_{10}.
 \end{equation}
Consequently, we can decompose the basic third-order invariants as follows:
 $$
R=R_{(0)}+R_{(2)}+R_{(4)}+R_{(6)},\quad
\nabla T=\nabla T_{(4)}+\nabla T_{(6)}+\nabla T_{(8)}+\nabla T_{(10)},
 $$
where the subscript $(k)$ indicates the weight, or equivalently the (multiple) submodule $V_k$ to which the
component belongs.
Indeed, the spaces where these tensors take values decompose into $\mathfrak{sl}(2)$-irreducibles as follows:
 $$
\mathfrak{gl}(2)\otimes\Lambda^2\tau^*=V_0\oplus 2V_2\oplus 2V_4\oplus V_6,\ \quad
\tau^*\otimes V_7=V_4\oplus V_6\oplus V_8\oplus V_{10}.
 $$
For a tensor $K$, the condition   $K\in V_l$ is equivalent to $C\cdot K=l(l+2)K$.
This gives a decomposition of $R$ and $\nabla T$ into eigenspaces of $C$.

By Proposition  5 of Sect. \ref{sec:lindeg}, system (\ref{main}) is linearly degenerate and integrable if and only if
the associated  Bryant connection is trivial: $T=R=0$ (this is also a consequence of Theorem \ref{decomp} below).
It turns out that in the general, not linearly degenerate case, the integrability is characterised by the condition that both $R$ and $\nabla T$ can be represented  as certain quadratic expressions in the torsion $T$.
To write down these expressions we introduce the following tensors
$T^2,\  T^2_\alpha,\  T^2_{\beta},\  T^2_{\gamma},\  T^2_{\delta}\in \tau\otimes\tau^*\otimes\Lambda^2\tau^*$
(the upper index 2 indicates that these tensors are quadratic in $T$, and  $\alpha, \beta, \gamma, \delta$
are labels, not indices):
 \begin{gather*}
(T^2)^k_{lij}=T^k_{la}T^a_{ij},\quad
(T^2_\alpha)^k_{lij}=T^k_{la}\Omega^{ab}T^c_{b[i}\Omega_{cj]},\quad
(T^2_\beta)^k_{lij}=T^k_{[ja}\Omega^{ab}T^c_{bl}\Omega_{ci]},\\
(T^2_\gamma)^k_{lij}=T^k_{[ia}\Omega^{ab}T^c_{bj]}\Omega_{cl},\quad
(T^2_\delta)^k_{lij}=\Omega^{ka}T^b_{al}\Omega_{bc}T^c_{ij};
 \end{gather*}
here square brackets  denote skew-symmetrization in  $i,j$, and $\Omega$ is the canonical almost symplectic structure
(note that these tensors are independent  of the choice of a conformal factor of $\Omega$).
Due to the identity $T^2=2\,T^2_\alpha$, the above formulae give 4 essentially different invariant tensors
that are quadratic in the second-order partial derivatives of  $f$ and $g$.
For every $\sigma\in\{\alpha,\beta,\gamma,\delta\}$ we decompose $T^2_\sigma$
into $\mathfrak{sl}(2)$-irreducibles, where apriori all weights $0\le l\le 10$, $l\in 2\mathbb{Z}$, are possible:
 $$
T^2_\sigma=T^2_\sigma{}_{(0)}+T^2_\sigma{}_{(2)}+T^2_\sigma{}_{(4)}+T^2_\sigma{}_{(6)}+T^2_\sigma{}_{(8)}
+T^2_\sigma{}_{(10)}.
 $$
We claim that in fact
 $$
T^2_\sigma{}_{(0)}=T^2_\sigma{}_{(4)}=T^2_\sigma{}_{(8)}=0\text{ for all }\sigma\in\{\alpha,\beta,\gamma,\delta\}.
 $$
Indeed, these  tensors are contractions of
$T\otimes T\in S^2V_7=V_2\oplus V_6\oplus V_{10}\oplus V_{14}$, which is then projected   to the right hand side of (\ref{dec-tten}).
Thus, $V_{14}$ disappears, and only components $V_k$ with $k= 2, 6,  10$ remain in the decomposition.

 \begin{theorem}
 \label{decomp}
For non-degenerate system (\ref{main}), the integrability  is equivalent to the following relations among the invariants of the associated
${\bf GL}(2,\mathbb{R})$ structure:
 \begin{gather*}
R_{(0)}=0,\ \ R_{(4)}=0,\ \ \nabla T_{(4)}=0,\ \ \nabla T_{(8)}=0,\ \ \nabla T_{(10)}=-28\,T^2_\alpha{}_{(10)},\\
R_{(2)}=\tfrac{44}3\,T^2_\alpha{}_{(2)}+2\,T^2_\beta{}_{(2)}-\tfrac{40}3\,T^2_\gamma{}_{(2)}-2\,T^2_\delta{}_{(2)},\\
R_{(6)}=-24\,T^2_\alpha{}_{(6)}-30\,T^2_\beta{}_{(6)}-60\,T^2_\gamma{}_{(6)}-24\,T^2_\delta{}_{(6)},\\
\nabla T_{(6)}=-8\,T^2_\alpha{}_{(6)}-8\,T^2_\beta{}_{(6)}-16\,T^2_\gamma{}_{(6)}-4\,T^2_\delta{}_{(6)}.
 \end{gather*}
 \end{theorem}

\noindent
Note that for $k=2,6$ there are 4 linearly independent components $T^2_\alpha{}_{(k)},T^2_\beta{}_{(k)},T^2_\gamma{}_{(k)},T^2_\delta{}_{(k)}$ of weight $k$,
and  $V_k$ has multiplicity 4 in the $\mathfrak{sl}(2)$-submodule
$\Pi_{(k)}=(\tau\otimes\tau^*\otimes\Lambda^2\tau^*)_{(k)}$.
Thus, it is expected that $R_{(k)}$ and $\nabla T_{(k)}$ decompose in this `basis'
(since the integrability conditions reduce third-order expressions to  second-order).
Similarly, for $k=10$ just one $T^2_\alpha{}_{(10)}$ is a basis of $\Pi_{(10)}$, and $\nabla T_{(10)}$
is expressed through it. However, the tensors $T^2_\sigma$ have no $V_k$-components for $k=0,4,8$, and for
 the corresponding third-order tensors we get simpler relations.

 \centerline{\bf Proof:}
 \medskip

The proof is computational. To verify these formulae we substitute  the integrability conditions (\ref{*}).  This leads to expressions that are quadratic in the
second-order partial derivatives of $f$ and $g$, with coefficients depending on their first-order derivatives (that is, functions on $J^1$). The computational complexity of
the output is quite high, however, the substitution of any `generic' point of $J^1$ readily gives zero (this verification
is done in Mathematica; we choose a rational point so that all computations are exact and rigorous).
Here  the choice of a generic point is irrelevant due to
the transitivity of ${\bf SL}(5)$-action on $J^1$, see discussion at the end of Sect. \ref{sec:equiv}.

Conversely, one has to verify that the above relations  imply {\it all} of the 40 integrability conditions (\ref{*}). Note that these  relations are linear  in third-order partial derivatives of $f$ and $g$, with
coefficients being functions on $J^1$. The rank of the matrix at third-order derivatives
is equal to 40 (one can again  restrict to any  generic point of $J^1$), and this implies the claim.

\medskip

\noindent {\bf Remark.} Some of the relations from Theorem \ref{decomp} hold identically for every fourfold
$X\subset {\bf Gr}(3,5)$, without using the integrability conditions. These include the  relations
 $$
R_{(0)}=0,\ T^2=2\,T^2_\alpha,\ T^2_\gamma{}_{(10)}=-T^2_\alpha{}_{(10)},\ T^2_\beta{}_{(10)}=T^2_\delta{}_{(10)}=0,
 $$
that play the role of obstructions to the embeddability of an abstract ${\bf GL}(2,\mathbb{R})$ structure into the Grassmannian ${\bf Gr}(3,5)$. It would be interesting to find a complete set of such obstructions.

\section{Concluding remarks}

This paper is a first step towards the general theory of integrability in Grassmann geometries. We gave a detailed characterisation of integrable systems $\Sigma(X)$ associated with fourfolds
 $X\subset {\bf Gr}(3, 5)$.

 \begin{itemize}

\item  We believe that some of our results can be generalised as follows.

\noindent (a) In the dimension $d=3$, the parameter space of non-degenerate integrable systems $\Sigma(X)$ associated with submanifolds  of codimension $n-3\geq 2$ in  $ {\bf Gr}(3, n)$ is finite-dimensional. Submanifolds $X$ corresponding to `generic' integrable systems
are not algebraic.

\noindent (b)  In higher dimensions $d\geq 4$, every non-degenerate integrable system $\Sigma(X)$ associated with a submanifold of codimension $n-d\geq 2$ in  $ {\bf Gr}(d, n)$ is necessarily linearly degenerate. Submanifolds $X$ corresponding to linearly degenerate integrable systems are  rational (generally, singular).

At the moment we are not aware of any $(n, d)$-independent approaches to the results of this kind.

 \item Some examples of integrable systems (\ref{main})  discussed in this paper can be interpreted as B\"acklund transformations:   on elimination of $v$, they lead to a second-order PDE for $u$, similarly, on elimination of $u$ they lead to a second-order PDE for $v$.  
 It would be interesting to obtain a classification of B\"acklund transformations.

 \item Our definition of integrability, based on the existence of holonomic trisecant submanifolds,  applies to any abstract ${\bf GL}(2, \bbbr)$ structure. It would be of interest to understand whether every integrable ${\bf GL}(2, \bbbr)$ structure is necessarily embeddable, that is, comes from a fourfold  $X\subset {\bf Gr}(3, 5)$. More generally, one may ask for a criterion of embeddability of an abstract ${\bf GL}(2, \bbbr)$ structure (see Remark at the end of Sect. \ref{sec:curvtor}).

 \item It is a true challenge to classify integrable systems that correspond to {\it algebraic} fourfolds $X\subset {\bf Gr}(3, 5)$. The homology class of any such $X$ can be represented as $a\sigma+b\eta$ where $a, b$ are nonnegative integers, and
 $\sigma, \eta$ are the standard four-dimensional Schubert cycles: $\sigma$ corresponds to 3-dimensional subspaces containing a fixed 1-dimensional subspace,  and $\eta$ corresponds to  3-dimensional subspaces that have 2-dimensional intersections with a fixed 3-dimensional subspace. Which values of $a$ and $ b$ are compatible with the requirement of integrability? The approach  of \cite{Bryant1, Robles} allows one to characterise algebraic $X$ (in some special homology classes) as integral manifolds of certain overdetermined exterior differential systems. Provided such characterisation is found for every $a$ and $ b$, it would be straightforward to intersect this differential system with our integrability conditions.

\end{itemize}

\section*{Acknowledgements}
We thank A Bolsinov, R Bryant, I Dolgachev, E Mezzetti, M Pavlov, A Prendergast-Smith and C Robles for clarifying discussions. We also thank the LMS for their support of BK to Loughborough making this collaboration possible.

\end{document}